\documentclass[a4paper,10 pt,reqno]{amsart}
\usepackage[latin1]{inputenc}
\usepackage[english]{babel}
\usepackage{amsmath, amsthm, amssymb, amsopn, amsfonts, amstext, stmaryrd, enumerate, color, mathtools, hyperref, mathrsfs, enumitem, url, cancel}
\usepackage[left=2.2cm,right=2.2cm,top=2.2cm,bottom=2.2cm]{geometry}
\numberwithin{equation}{section}
\usepackage{dsfont}

\hypersetup{
    colorlinks=true,
    linkcolor=blue!75!black,
    citecolor=red!50!black,
    urlcolor=green!50!black,
    linktoc=all
}

\usepackage{colortbl}
\usepackage{tikz}

\newcommand{\C}{\mathcal{C}}

\newcommand{\I}{\mathcal{I}}

\newcommand{\R}{\mathbb{R}}
\newcommand{\RN}{\mathbb{R}^N}

\newcommand{\dist}{{\mbox{\normalfont dist}}}

\newcommand{\Lfrac}{(-\Delta)^{s}}
\newcommand{\bu}{\mathbf{u}}

\newcommand{\vertiii}[1]{{\left\vert\kern-0.25ex\left\vert\kern-0.25ex\left\vert #1
    \right\vert\kern-0.25ex\right\vert\kern-0.25ex\right\vert}}

\newcommand{\re}{\mathbb{R}}
\newcommand{\ren}{\re^N}

\newcommand{\io}{\int\limits_\O}

\renewcommand{\a }{\alpha }
\renewcommand{\b }{\beta }
\renewcommand{\d }{\delta }
\newcommand{\D }{\Delta }

\newcommand{\n }{\nabla }

\newcommand{\s }{\sigma }

\renewcommand{\O }{\Omega }

\renewcommand{\epsilon} {\varepsilon}

\DeclareMathOperator{\diam}{diam}

\def\XXint#1#2#3{{\setbox0=\hbox{$#1{#2#3}{\int}$ }
\vcenter{\hbox{$#2#3$ }}\kern-.6\wd0}}

\theoremstyle{plain}
\newtheorem{theorem}{Theorem}[section]
\newtheorem{proposition}[theorem]{Proposition}
\newtheorem{lemma}[theorem]{Lemma}
\newtheorem{corollary}[theorem]{Corollary}

\theoremstyle{definition}
\newtheorem{definition}{Definition}

\theoremstyle{remark}
\newtheorem{remark}[theorem]{Remark}

\renewcommand{\le}{\leqslant}
\renewcommand{\leq}{\leqslant}

\renewcommand{\geq}{\geqslant}

\begin{document}

\title[Global fractional Calder\'on--Zygmund type regularity]{Global fractional Calder\'on--Zygmund type regularity}

\author[B. Abdellaoui, A. J. Fern\'andez, T. Leonori and A. Younes]{Boumediene Abdellaoui, Antonio J. Fern\'andez, Tommaso Leonori and Abdelbadie Younes}


\address{
\newline
\textbf{{\small Boumediene Abdellaoui}}
\vspace{0.1cm}
\newline \indent D\'epartement de Math\'ematiques, Universit\'e Abou Bakr Belka\"id, Tlemcen, 13000, Algeria}
\email{boumediene.abdellaoui@inv.uam.es}

\address{
\vspace{-0.25cm}
\newline
\textbf{{\small Antonio J. Fern\'andez}}
\vspace{0.1cm}
\newline \indent Instituto de Ciencias Matem\'aticas, Consejo Superior de Investiaciones Cient\'ificas, 28049 Madrid, Spain
\newline \indent Departamento de Matem\'aticas, Universidad Aut\'onoma de Madrid, 28049, Madrid, Spain}
\email{antonioj.fernandez@uam.es}

\address{
\vspace{-0.25cm}
\newline
\textbf{{\small Tommaso Leonori}}
\vspace{0.1cm}
\newline \indent Dipartimento di Scienze di Base e Applicate per l'Ingegneria, Sapienza Universit\`a di Roma, 00161 Roma, Italy}
\email{tommaso.leonori@uniroma1.it}

\address{
\vspace{-0.25cm}
\newline
\textbf{{\small Abdelbadie Younes}}
\vspace{0.1cm}
\newline \indent D\'epartement de Math\'ematiques, Universit\'e Abou Bakr Belka\"id, Tlemcen, 13000, Algeria}
\email{abdelbadieyounes@gmail.com}

\begin{abstract}
We obtain a global fractional Calder\'on-Zygmund type regularity theory for the fractional Poisson problem. More precisely, for $\Omega \subset \RN$, $N \geq 2$, a bounded domain with boundary $\partial \Omega$ of class $C^2$, $s \in (0,1)$ and $f \in L^m(\Omega)$ for some $m \geq 1$, we consider the problem
$$
\left.
\begin{aligned}
(-\Delta)^s u = f \quad \textup{in } \Omega,  \qquad\ u = 0 \quad \textup{in }   \RN \setminus \Omega,
\end{aligned}
\right.
$$
and, according to $m$, we find the values of $s \leq t < \min\{1,2s\}$ and of $1 < p < +\infty$ such that $u \in L^{t,p}(\RN)$ and such that $u \in W^{t,p}(\RN)$.
\medbreak
\smallbreak
\end{abstract}

\maketitle
$ $
\vspace{-1.17cm}

\section{Introduction} \label{introduction}
\noindent  The aim of the present paper is to obtain a global fractional Calder\'on-Zygmund type regularity theory for the fractional Poisson  problem  with  homogeneous  Dirichlet boundary conditions, namely
\begin{equation} \tag{P} \label{main0}
\left\{
\begin{aligned}
\Lfrac u & = f, &&  \textup{ in } \Omega,\\
u & = 0, && \textup{ in }   \RN \setminus \Omega.
\end{aligned}
\right.
\end{equation}
We assume here and in the rest of the paper that $\Omega \subset \RN$, $N\geq 2$, is a bounded domain with boundary $\partial \Omega$ of class $C^2$ and that $f \in L^m(\Omega)$ for some $m \geq 1$.  Also, for $s \in (0,1)$, we denote by $(-\Delta)^s$ the fractional Laplacian which, for all $u \in C_c^{\infty}(\RN)$, is defined by
\[ (-\Delta)^s u(x):= a_{N,s}\textup{ p.v.}\int_{\RN} \frac{u(x)-u(y)}{|x-y|^{N+2s}} dy = a_{N,s} \lim_{\epsilon \to 0^{+}} \int_{\RN \setminus B_{\epsilon}(x)} \frac{u(x)-u(y)}{|x-y|^{N+2s}} dy,\]
where $a_{N,s} := 2^{2s} s \pi^{-\frac{N}{2}} \frac{ \Gamma \left( \frac{N}{2}+s \right)}{ \Gamma(1-s) }$ is a normalization constant and $\Gamma$ denotes the Euler's gamma function.

Let us  rephrase our  main goal as follows: we want to find the fractional regularity of the (unique) solution to \eqref{main0} according to the values of $m$.
In particular, we focus our attention on the regularity in suitable function spaces such as  the Bessel potential spaces $ L^{t,p}(\RN)$ and the fractional Sobolev spaces $W^{t,p}(\RN)$ with  $s \leq t < \min\{1,2s\}$ and  $1 < p < + \infty$ (see \textit{Function spaces} below for a reminder of their definitions).

Before stating and describing our main results, we want to illustrate what kind of regularity can be expected through an example.
Let us consider the \textit{fractional torsion problem}
\begin{equation} \tag{T} \label{torsion-problem}
\left\{
\begin{aligned}
\Lfrac u & = 1, &&  \textup{ in } B_1(0),\\
u & = 0, && \textup{ in } \RN \setminus B_1(0).
\end{aligned}
\right.
\end{equation}
It is well known (see e.g. \cite[Theorem 1]{Dyda}) that the unique solution to \eqref{torsion-problem} is
\begin{equation} \label{torsion-function}
\mathbf{u}(x) :=  \left\{
\begin{aligned}
& \ \frac{2^{-2s}\,\Gamma\big(\frac{N}{2}\big)}{\Gamma\big(\frac{N}{2} + s \big)  \Gamma(1+s)} (1-|x|^2)^s, && \textup{ in } B_1(0), \\
& \ 0, && \textup{ in } \RN \setminus B_1(0).
\end{aligned}
\right.
\end{equation}
Having at hand this explicit expression, it is easy to verify that $\bu \in W^{1,p}(\RN)$ if and only if $1 \leq p < \frac{1}{1-s}$. This shows that, even if the right hand side $f$ (as well as the boundary of the domain $\Omega$) in the fractional Poisson problem \eqref{main0} is smooth, one cannot expect in general the solution  to be smooth up to the boundary. This is in striking contrast with the regularity results for the Poisson problem  associated to the Laplacian operator (see for example the classical book by D. Gilbarg and N. S. Trudinger \cite{G_T_2001_S_Ed}). Regardless, such peculiarity of the fractional Laplacian does not exclude the possibility of having a global Calder\'on--Zygmund type regularity theory in fractional Sobolev spaces and Bessel potential spaces. For instance, using the explicit expression of $\bu$ given in \eqref{torsion-function} and a result by B. Dyda (see \cite[Theorem 1]{Dyda}), it is easy to prove the following result:

\begin{proposition} \label{integrability-t-Laplacian-intro} Let $\mathbf{u}$ be the unique solution to \eqref{torsion-problem}. Then:
\begin{enumerate}
\item[$i)$] If $0 < t < s$, then $(-\Delta)^{\frac{t}{2}} \bu \in L^{p}(\RN)$ for all $1 \leq p \leq \infty$.
\item[$ii)$]  $(-\Delta)^{\frac{s}{2}} \bu \in L^p(\RN)$ for all $1 \leq p < \infty$, but  $(-\Delta)^{\frac{s}{2}}\bu \not \in L^{\infty}(B_1(0))$.
\item[$iii)$] If $s < t < \min\{1,2s\}$, then $(-\Delta)^{\frac{t}{2}} \bu \in L^p(\RN)$ for all $1 \leq p < \frac{1}{t-s}$, but $(-\Delta)^{\frac{t}{2}}\bu \not \in L^p(B_1(0))$ for any $ \frac{1}{t-s} \leq p \leq \infty$.
\end{enumerate}
\end{proposition}

The analysis of the solution  to  the \textit{fractional torsion problem} \eqref{torsion-problem} deeply relies on the explicit expressions of $\bu$ and $(-\Delta)^{\frac{t}{2}}\bu$. However, in general, we do not explicitly know the expression neither of the solution to \eqref{main0} nor of its $\frac{t}{2}$--Laplacian. Keeping in mind Proposition \ref{integrability-t-Laplacian-intro}, we aim to derive a  global fractional Calder\'on--Zygmund type regularity theory  for more general domains $\Omega$ and data functions $f$, where no explicit expression of the solution is known.   We first deal with $f$ ``very integrable'' and prove the following sharp result.

\begin{theorem} \label{calderon-zygmund-very-integrable-intro}
Let $s \in (0,1)$, $s \leq t < \min\{1,2s\}$ and let $u$ be the (unique) solution to \eqref{main0} with $f \in L^{m}(\Omega)$ for some $m > \frac{N}{2s-t}$. Then: \smallbreak
\begin{itemize}
\item[i)] For all $1 \leq p < \infty$, there exists $C > 0$ such that
$$
\|(-\Delta)^{\frac{s}{2}}u\|_{L^p(\RN)} \leq C \|f\|_{L^{m}(\Omega)}.
$$
\item[ii)] For all $1 \leq p < \frac{1}{t-s}$, there exists $C > 0$ such that
$$
\|(-\Delta)^{\frac{t}{2}}u\|_{L^p(\RN)} \leq C \|f\|_{L^{m}(\Omega)}.
$$
\end{itemize}
Here, $C > 0$ are constants depending only on $N$, $s$, $t$, $p$, $m$ and $\Omega$.
\end{theorem}

\begin{remark} $ $
\begin{itemize} \item[a)] Proposition \ref{integrability-t-Laplacian-intro} guarantees that the previous result (see also Theorem \ref{calderon-zygmund-very-integrable}) is sharp.
\item[b)] In the case where $t = s$ we use the convention $\frac{1}{t-s} = +\infty$. Hence, in that situation, $i)$ and $ii)$ of the previous result coincide.
\item[c)] The previous result is a particular case of the more general Theorem \ref{calderon-zygmund-very-integrable}, where  we also analyze the (somehow simpler) case where $0 < t < s$.
\end{itemize}
\end{remark}

The analysis of the global Calder\'on--Zygmund type regularity we perform goes beyond the case of ``regular'' right hand side. Next, we pursue this analysis with a right hand side with ``low summability''. Let us introduce some notation. For $s \in (0,1)$, $s \leq t < \min\{1,2s\}$, $\epsilon \in (0,1)$ arbitrary small and $1 \leq m < \frac{N}{2s-t}$, we define
\begin{equation*}
m_{\star}(s,t):=  \max \Big\{ 1, \, \frac{N}{N+s-N(t-s)}+\epsilon \Big\},
\end{equation*}
and
\begin{equation*}
p^{\star}(m,s,t):= \left\{
\begin{aligned}
& \frac{mN}{N-ms},  && \textup{ for } t = s,\\
& \min\Big\{ \frac{mN}{N-ms+mN(t-s)}, \frac{1}{t-s} \Big\}, \quad && \textup{ for } t > s.
\end{aligned}
\right.
\end{equation*}

\bigbreak
\noindent Our main global Calder\'on--Zygmund type regularity result concerning ``low summability'' $f$ reads as follows:

\begin{theorem} \label{calderon-zygmund-low-integrable}
Let $s \in (0,1)$, $s \leq t < \min\{1,2s\}$ and let $u$ be the (unique) solution to \eqref{main0} with $f \in L^m(\Omega)$ for some $m_{\star}(s,t) \leq m < \frac{N}{2s-t}$. Then, for all $1 \leq p < p^{\star}(m,s,t)$, there exists $C > 0$ (depending only on $N$, $s$, $t$, $p$, $m$ and $\Omega$) such that
$$
\|(-\Delta)^{\frac{t}{2}}u\|_{L^p(\RN)} \leq C\, \|f\|_{L^m(\Omega)}\,.
$$
\end{theorem}

As a byproduct of our approach to prove Theorems \ref{calderon-zygmund-very-integrable-intro} and \ref{calderon-zygmund-low-integrable}, we also estimate the behaviour of the $\frac{t}{2}$--Laplacian of the solution to \eqref{main0} at  the boundary of the domain $\Omega$.

\begin{theorem} \label{weighted-distance-intro}
Let $s \in (0,1)$, $s \leq t < \min\{1,2s\}$ and let $u$ be the (unique) solution to \eqref{main0} with $f \in L^{m}(\Omega)$ for some $m \geq 1$.
\smallbreak
\begin{itemize}
\item[i)] If $m >\frac{N}{2s-t}$, there exists $C > 0$ such that
$$
\big\| \, |\log \delta|^{-1}\, (-\Delta)^{\frac{s}{2}} u \big\|_{L^{\infty}(\Omega)} + (t-s) \big\| \delta^{t-s} (-\Delta)^{\frac{t}{2}}u\big\|_{L^{\infty}(\Omega)}  \leq C\, \|f\|_{L^m(\Omega)}.
$$
\smallbreak
\item[ii)] If $1\leq m <\frac{N}{2s-t}$, for all $1 \leq p < \frac{mN}{N-ms}$ and $1 \leq q < \frac{mN}{N-m(2s-t)}$, there exists $C > 0$ such that
$$
\big\| \, |\log \delta|^{-1}\, (-\Delta)^{\frac{s}{2}} u \big\|_{L^{p}(\Omega)} + (t-s) \big\| \delta^{t-s} (-\Delta)^{\frac{t}{2}} u \big\|_{L^q(\Omega)} \leq C\, \|f\|_{L^m(\Omega)}.
$$
\end{itemize}
Here, $\delta(x) := \dist(x,\partial \Omega)$ and  $C > 0$ are constants depending only on $N$, $s$, $t$, $p$, $q$, $m$ and $\Omega$.
\end{theorem}

\bigbreak

\begin{remark} $ $
\begin{itemize}
\item[a)]
A straightforward outcome of the above result is the following local version of the Calder\'on--Zygmund type regularity theory.  For any $\omega \subset \subset \O$, it follows that:

\smallbreak
\begin{itemize}
\item[$i)$]  If $m > \frac{N}{2s-t}$, there exists $C_{\omega} > 0$ such that
$$
\big\|   (-\Delta)^{\frac{s}{2}} u \big\|_{L^{\infty}(\omega)} + (t-s)\big\|   (-\Delta)^{\frac{t}{2}}u\big\|_{L^{\infty}(\omega)}  \leq C_\omega \, \|f\|_{L^m(\Omega)}.
$$
\smallbreak
\item[$ii)$] If $1\leq m < \frac{N}{2s-t}$, for all $1 \leq p < \frac{mN}{N-ms}$ and $1 \leq q < \frac{mN}{N-m(2s-t)}$, there exists $C_{\omega} > 0$ such that
$$
\big\|   (-\Delta)^{\frac{s}{2}} u \big\|_{L^{p}(\omega)} + (t-s)  \big\|   (-\Delta)^{\frac{t}{2}}u\big\|_{L^{q}(\omega)} \leq C_\omega\, \|f\|_{L^m(\Omega)}.
$$
\end{itemize}
\item[b)] The above theorem  is a particular case of the slightly more general Theorem \ref{prop-distance-weight}.
\end{itemize}
\end{remark}

\bigbreak
As a consequence of Theorems \ref{calderon-zygmund-very-integrable-intro} and \ref{calderon-zygmund-low-integrable}, we get the following regularity in fractional Sobolev spaces.

\begin{corollary} \label{cor-Wsp}
Let $s \in (0,1)$ and let $u$ be the unique solution to \eqref{main0} with $f \in L^m(\Omega)$. \smallbreak
\begin{itemize}
\item[i)] If $1 \leq m < \frac{N}{s}$, then, for all $1 < p < \frac{mN}{N-ms}$, there exists $C > 0$ such that
$$
\|u\|_{W^{s,p}(\RN)} \leq C \|f\|_{L^m(\Omega)}.
$$
\item[ii)] If $m > \frac{N}{s}$, then, for all $1 < p < +\infty$, there exists $C > 0$ such that
$$
\|u\|_{W^{s,p}(\RN)} \leq C \|f\|_{L^m(\Omega)}.
$$
\end{itemize}
Here, $C > 0$ are positive constants depending only on $N$, $s$, $p$, $m$ and $\Omega$.
\end{corollary}

\begin{remark}The above result is a particular case of the more general results Corollary \ref{cor2-high-integrability} and Corollary \ref{cor2-low-integrability}.
\end{remark}

Now that our main regularity results are stated, let us locate our work with respect to the existing literature. To that end, let us denote by $u$ the (unique) solution to \eqref{main0}. Regarding the regularity of $u$, there exist mainly two different subjects: \textit{local (or interior)} regularity and \textit{global (or up to the boundary)} regularity.

Let us first focus in the existing literature concerning \textit{local} regularity. This direction is nowadays rather well understood so we just mention a few works dealing with different issues. On one hand, let us mention that the $W^{2s-\epsilon,2}_{loc}(\Omega)$ and $L^{2s,p}_{loc}(\Omega)$ regularity of $u$ have been respectively analyzed in \cite{cozzi} (where the author deals with more general operators) and \cite{BWZ}. On the other hand, in the very recent paper \cite{MSY}, the \textit{local} summability of the fractional derivatives of $u$ of order $t$, with $s\leq t < \min\{ 1,2s\}$, has been treated. Note that in \cite{MSY} the authors deal with general operators, whose kernel is a $C^{\alpha}$--perturbation of the fractional Laplacian one. Even more general operators have been now treated in \cite{N22}. We should also mention that related results concerning \textit{local} Sobolev type regularity have been proved in \cite{N20, SSS15, Y17}. Finally, we refer to \cite[Proposition 4.3]{DPSV}, where sharp  Schauder type estimates have been proved. Let us also mention the recent work \cite{fall} for the analysis of more general operators.

As already observed, the \textit{global} regularity is much more involved. Let us first recall that the {\textit{global}} $L^p$--regularity of $u$ has been completely understood, even for more general operators than the fractional Laplacian (see for instance \cite{LPPS2015}). Moreover, when the datum  $f$ is ``regular'' (say $f \in L^m(\Omega)$ with  $m > \frac{N}{2s}$), sharp {\textit{global}} H\"older--type regularity results have been obtained. We refer to \cite{AG23, GRUBB, RO-S14-2, Sil} for fundamental contributions in this direction.
Let us  in addition  mention that the behavior of $u$ near the boundary of the domain has also been studied. We refer to \cite{GRUBB, RO-S14-2} for the analysis of the behaviour of $ \delta^{-s} u$ and to \cite{F-J-2021} for the analysis of the behaviour of $\delta^{1-s} \nabla u$. Despite the above mentioned works, very few results dealing with the \textit{global} summability of the fractional derivatives of $u$ or with their behaviour near the boundary of the domain are available in the literature.
Among these, let us point out \cite{KMS}, where the authors established {\textit{global}} summability results for the fractional derivatives of order $\sigma$ with $0 < \sigma < s$. Note that they deal with operators much more general than the fractional Laplacian and that their proofs rely on the use of the Wolff potential.
The same type of results, using a completely different approach, has been proved in  \cite{AbdAttBen}. Concerning \textit{global} fractional Sobolev regularity, we should also mention \cite[Section 2]{BLN21} and \cite[Section 2 and 3]{BNS19}, where \textit{global} regularity results in weighted Sobolev spaces have been proved, and \cite{BN23}, where the authors deal with $L^2$ right hand sides and obtain \textit{global} regularity results in Besov spaces.

The main novelty of our results with respect to the existing literature is that we are able to understand the \textit{global} summability of the fractional derivatives of $u$ for any order $t$, with $s\leq t <  \min\{ 1,2s\}$, and their behaviour near the boundary of the domain. On one hand, for the ``particular case'' of the fractional Laplacian, our results extend to \textit{global} the recent \textit{local} Calder\'on--Zygmund type regularity theory established in \cite{MSY, N22}. On the other hand, we obtain sharp bounds for $\delta^{t-s} (-\Delta)^{\frac{t}{2}}u$ near the boundary of the domain. Our analysis of $\delta^{t-s} (-\Delta)^{\frac{t}{2}}u$ goes on the line of the works \cite{F-J-2021,GRUBB, RO-S14-2} but our results are new.

Next, we give some ideas of the proofs of Theorems \ref{calderon-zygmund-very-integrable-intro} and \ref{calderon-zygmund-low-integrable}. Let us first specify the notion of solution to \eqref{main0} that we use in this work. We denote by $G_s : \R_{\ast}^{2N} \to \R$ the Green kernel associated to $(-\Delta)^{s}$ in $\Omega$ and introduce the \textit{solution map}
$$
\mathbb{G}_s : L^1(\Omega) \to M^{\frac{N}{N-s}}(\Omega), \qquad f \mapsto \mathbb{G}_s[f] := \int_{\Omega} G_s(x,y) f(y)\, dy,
$$
where $M^{p}(\Omega)$ with $p\geq1$ denotes the Marcinkiewicz space (see \eqref{marc} for the definition).  Note that the \textit{solution map} is well-defined thanks to the existence and uniqueness of a solution to \eqref{main0} with $L^1$--data and by its representation formula (see e.g. \cite{C_V_2014-JDE}).

\begin{definition} 
We say that $u \in L^1(\Omega)$ is a (very weak) solution to \eqref{main0} if
\begin{equation} \label{representation-formula}
u(x) = \mathbb{G}_s[f](x) = \int_{\Omega} G_s(x,y) f(y) dy, \quad \textup{ for  a.e. } x \in \R^N.
\end{equation}
\end{definition}

\bigbreak
\begin{remark} $ $
\begin{itemize}
\item[a)] We will sometimes use the notation $u:= \mathbb{G}_s[f]$ to indicate that $u$ is the solution to \eqref{main0}.
\smallbreak
\item[b)] Let us define
\begin{equation*} 
\mathbb{X}^s(\Omega) := \Big\{\phi\in C^s(\ren)\,: \, \phi (x) = 0 \textup{ for all }x \in \RN \setminus \Omega \textup{ and } (-\Delta)^s \phi \in L^{\infty}(\Omega)\Big\}.
\end{equation*}
It is more common in the literature to use the following notion of (very weak) solution: $u \in L^1(\Omega)$ is a \textit{very weak solution} to \eqref{main0} if $u \equiv 0$ in $\RN \setminus \Omega$ and
$$
\io u(-\Delta )^s\phi\, dx =\io \phi(x) f(x) dx, \quad \textup{ for all }  \phi  \in \mathbb{X}^s(\Omega)\,.
$$
However, as explained for instance in \cite[Remark 2.11]{AP-KPZ-2018}, both notions are equivalent. If one prefers so, one could identify our definition of solution as the representation formula for \textit{very weak solutions}. \smallbreak 
\item[c)] For all $y \in \Omega$, let $\delta_y$ be the Dirac mass centred at $y$. It is well-known that the Green function $G_s (\cdot,y)$ is,  at least formally, a solution to
\begin{equation*}
\left\{
\begin{aligned}
(-\Delta)^s_x G_s (x,y) & = \delta_y, &&  \textup{ in } \Omega,\\
G_s (x,y) & = 0, && \textup{ in }   \ren\setminus \Omega \,.
\end{aligned}
\right.
\end{equation*}
\end{itemize}
\end{remark}

\bigbreak
Having at hand the notion of solution, we are now ready to describe the main ideas in the proofs of Theorems \ref{calderon-zygmund-very-integrable-intro} and \ref{calderon-zygmund-low-integrable}. Let $s \in (0,1)$, $s \leq t < \min\{1,2s\}$ and let $u$ be the (unique) solution to \eqref{main0} with $f \in L^m(\Omega)$ for some $m \geq m_{\star}(s,t)$. As a first step, we notice that (see Remark \ref{f1.6} for more details) the representation formula \eqref{representation-formula} provides a suitable representation formula for the $\frac{t}{2}$-Laplacian of $u$:
\begin{equation} \label{representation-formula-t-half-laplacian}
(-\Delta)^{\frac{t}{2}}u(x) = \int_{\Omega} (-\Delta)^{\frac{t}{2}}_x G_s(x,y) f(y) dy, \quad \textup{ for a.e. } x \in \Omega.
\end{equation}

\noindent We are then lead to analyze the pointwise behaviour of the $\frac{t}{2}$--Laplacian of the Green kernel $G_s$. Our main result in that direction reads as follows.

\begin{theorem} \label{main-th-Green}
Let $s \in (0,1)$ and $s \leq t < \min\{1,2s\}$. There exists a constant $C > 0$ (depending only on $N,s,t$ and $\Omega$) such that
\begin{equation}\label{grest}
\big|(-\Delta)_x^{\frac{t}{2}} G_s(x,y)\big| \leq \frac{C}{|x-y|^{N-(2s-t)}} \left( \big|\log |x-y|\big| + |\log\delta(x)| + \frac{|x-y|^{t-s}}{\d^{t-s} (x)} \right), \quad \textup{ for a.e. } x, y \in \Omega.
\end{equation}
Here, $\delta(x):= \dist(x,\partial\Omega)$.
\end{theorem}

Once we have at hand the representation formula for the $\frac{t}{2}$-Laplacian \eqref{representation-formula-t-half-laplacian} and Theorem \ref{main-th-Green}, the proofs of Theorems \ref{calderon-zygmund-very-integrable-intro} and \ref{calderon-zygmund-low-integrable} are actually standard. Indeed, combining these two elements, we will obtain a pointwise estimate for the $\frac{t}{2}$-Laplacian of $u$ (see Lemma \ref{declem} below) and then Theorems \ref{calderon-zygmund-very-integrable-intro} and \ref{calderon-zygmund-low-integrable} will follow from classical tools of harmonic analysis. We refer to Sect. \ref{sec5} for more details. Let us also mention that the strategy to prove Theorem  \ref{weighted-distance-intro} follows the same lines. Certainly, looking carefully to the pointwise estimate established in Theorem \ref{main-th-Green}, one can understand the role played by the distance function in Theorem  \ref{weighted-distance-intro}.

In a slightly different direction, we also prove pointwise estimates for the \textit{Riesz $t$-gradient} of the Green kernel $G_s$. Let us recall that, for all $\phi \in C_c^{\infty}(\RN)$, the \textit{Riesz fractional gradient of order $t$} is given by
\begin{equation*}
\n^t \phi(x):= \mu_{N,t} \int_{\RN} \frac{(x-y)(\phi(x)-\phi(y))}{|x-y|^{N+t+1}}\, dy,
\end{equation*}
where $ \mu_{N,t} := 2^t \pi^{-\frac{N}{2}} \frac{\Gamma (\frac{N+t+1}{2} )}{\Gamma ( \frac{1-t}{2} )}$ is a normalization constant. Modifying appropriately the proof of Theorem \ref{main-th-Green}, we obtain the following result.

\begin{theorem} \label{gradGreen-intro}
Let $s \in (0,1)$ and $s \leq t < \min\{1,2s\}$. There exists a constant $C > 0$ (depending only on $N$, $s$, $t$ and $\Omega$) such that
\begin{equation*}
\big|\nabla_x^{t} G_s(x,y)\big| \leq \frac{C}{|x-y|^{N-(2s-t)}} \left( \big|\log |x-y|\big| + |\log\delta(x)| + \frac{|x-y|^{t-s}}{\d^{t-s} (x)} \right), \quad \textup{ for a.e. } x, y \in \Omega.
\end{equation*}
Here, $\delta(x):= \dist(x,\partial\Omega)$.
\end{theorem}

Having at hand Theorem \ref{gradGreen-intro}, one can argue as in the proofs of Theorems \ref{calderon-zygmund-very-integrable-intro} and \ref{calderon-zygmund-low-integrable} and obtain their counterpart for the \textit{Riesz fractional gradient of order $t$}. We refer to Theorems \ref{CZ-gradient1} and \ref{CZ-gradient2} for the precise statements of these results.

Let us provide as well some ideas of the proofs of Theorems \ref{main-th-Green} and \ref{gradGreen-intro}. Both proofs mainly rely on two ingredients. On one hand, we use some classical   estimates on the Green kernel $G_s$ and on its gradient (see Lemma \ref{Green-estimates-lemma}). On the other hand, we use some new estimates for singular integrals with suitable weights that involve the distance function $\delta$ (see Lemma \ref{lemma-Tobias}). Indeed, the proof of Theorem \ref{main-th-Green} follows combining these two ingredients with a suitable decomposition for the $\frac{t}{2}$-Laplacian (and  the $t$-gradient, as far as Theorem \ref{gradGreen-intro} is concerned) of the Green kernel $G_s$.


As an application of our \textit{global} Calder\'on-Zygmund type regularity theory,  in our companion paper  \cite{AFLY-2021}, we get new existence results for nonlinear fractional Laplacian problems with nonlocal ``gradient terms''. In other words, we prove new existence results for deterministic stationary KPZ--type equations involving nonlocal-nonlinear ``gradient'' terms. More precisely, we analyze equations of the form
\begin{equation*}
\left\{
\begin{aligned}
(-\Delta)^s u & = \mu(x) \, |\mathbb{D} (u)|^q + \lambda h(x),  \quad  & \mbox{ in } \Omega\\
u & =0,  && \mbox{ in } \RN \setminus \Omega,
\end{aligned}
\right.
\end{equation*}
where $q > 1$ and $\lambda > 0$ are real parameters, $h$ belongs to a suitable Lebesgue space, $\mu \in L^{\infty}(\Omega)$ and $\mathbb{D}$ represents different nonlocal ``gradient'' terms. These  results provide (almost) complete answers to some of the open problems posed in \cite[Section 6]{AP-KPZ-2018} and \cite[Section 7]{AF-2020}.
 
The paper is organized as follows. In the rest of the introduction, we gather some notation and recall the definition of the function spaces used in this work. Sect. \ref{estimates-section} is devoted to prove several integral estimates that will be key in the rest of the paper. We also state there some tools that will be useful in the subsequent sections. Having at hand the integral estimates of Sect. \ref{estimates-section}, we prove Theorems \ref{main-th-Green} and \ref{gradGreen-intro} in Sect. \ref{sec4}. In Sect. \ref{sec5}, we use Theorem \ref{main-th-Green} to establish a \textit{global} Calder\'on-Zygmund type regularity theory for the fractional Poisson problem. In particular, we prove Theorems \ref{calderon-zygmund-very-integrable-intro}, \ref{calderon-zygmund-low-integrable} and \ref{weighted-distance-intro} and Corollary \ref{cor-Wsp}. The paper is closed by an appendix where we prove the elementary Proposition \ref{integrability-t-Laplacian-intro} and derive some immediate corollaries.


\subsection*{Notation.} We gather here some notation that will be used throughout the paper. For a set $\omega \subset \RN$, we define $\omega^c:= \RN \setminus \omega$. Also, we denote by $\d$ the distance to $\partial \Omega$, i.e. $\delta(x) = \dist(x,\partial \Omega)$,
and define
\begin{equation}\label{erre}
R:=  \frac13 + \frac43 \big(\diam(\Omega) + \dist(0,\Omega)\big)
\qquad \mbox{  so that } \quad \Omega \subset\subset B_{R}(0)\,.
\end{equation}
Note that, with such a choice of $R$, we have the following inequality that we use widely  in the sequel
\begin{equation}\label{erre2}
|x-y| \geq \frac14 \big(1 + |x|\big) \quad \textup{ for all } y \in \Omega \textup{ and all } x \in B_R^c(0) \,.
\end{equation}
Since it will play a role in the proof of Theorem \ref{main-th-Green}, we emphasize that, by $(-\Delta)^{\gamma}_{\Omega}$ we denote the regional fractional Laplacian which, for all $u \in C_c^{\infty}(\RN)$, is defined by
$$
(-\Delta)_{\Omega}^{\gamma} u (x) := a_{N,\gamma}\, {\rm{p.v.}} \int_{\Omega} \frac{u (x) - u (z)}{|x-z|^{N+2\gamma}}\, dz.
$$
 
It is also worth to emphasize here that, in the rest of the paper, we will omit the normalization constant $a_{N,\gamma}$ and the principal value sense (``p.v.'') appearing in the definitions of $(-\Delta)^\gamma$ and $(-\Delta)_{\Omega}^\gamma$. We will omit as well the normalization constant $\mu_{N,\gamma}$ appearing in the definition of $\nabla^{\gamma}$.
 
Finally, let us define $\R^{2N}_{\ast} := \{ (x,y) \in \R^{2N} : x \neq y\}$ and stress that, by $\sigma(\mathbb{S}^{N-1})$, we denote the $(N-1)$-dimensional measure of the unit sphere in $\RN$.

\subsection*{Function spaces.}
We collect here the definitions of the functional spaces involved in our results. First of all, recall that, for all $s \in (0,1)$ and $ 1 \leq p < +\infty$, the fractional Sobolev space $W^{s,p}(\RN)$ is defined as
$$
W^{s,p}(\RN):= \left\{ u \in L^p(\RN) : \iint_{\R^{2N}} \frac{|u(x)-u(y)|^p}{|x-y|^{N+sp}} dx dy < +\infty \right\}.
$$
It is a Banach space endowed with the usual norm
$$
\|u\|_{W^{s,p}(\RN)} = \left( \|u\|_{L^p(\RN)}^p + \iint_{\R^{2N}} \frac{|u(x)-u(y)|^p}{|x-y|^{N+sp}} dx dy \right)^{\frac{1}{p}}.
$$

\noindent Next, we remind that, for any $s \in (0,1)$ and $1 \leq p < +\infty$, the Bessel potential space  is defined as
$$
L^{s,p}(\RN) \ :=   \ \overline{ \big\{u \in C_c^{\infty}(\RN)\big\} }^{\, \vertiii{\cdot}_{L^{s,p}(\RN)}},
$$
where
$$
\vertiii{u}_{L^{s,p}(\RN)} = \|(1-\Delta)^{\frac{s}{2}}u\|_{L^p(\RN)} \quad \textup{ and } \quad (1-\Delta)^{\frac{s}{2}}u = \mathcal{F}^{-1} ( (1+|\cdot|^2)^{\frac{s}{2}} \mathcal{F} u), \quad \forall\ u \in C_c^{\infty}(\RN).
$$
Let us stress that, in the case where $s \in (0,1)$ and $1 < p < +\infty$,
$$
\|u\|_{L^{s,p}(\RN)} := \|u\|_{L^p(\RN)} + \|(-\Delta)^{\frac{s}{2}}u\|_{L^p(\RN)}
$$
is an equivalent norm for $L^{s,p}(\RN)$ (see e.g. \cite[page 5]{AC2020} for a precise explanation of this fact). Let us as well recall that, for all $0 < \epsilon < s < 1$ and all $1 < p < +\infty$, by \cite[Theorem 7.63, (g)]{adams}, we have
$$
L^{s+\epsilon,p}(\RN) \subset W^{s,p}(\RN) \subset L^{s-\epsilon,p}(\RN).
$$
Finally, we emphasize that, for $1 \leq p < +\infty$, $M^{p}(\Omega)$ denotes the Marcinkiewicz space
(or weak--$L^{p}$ space). We recall here that
\begin{equation} \label{marc}
[f]_{M^p(\Omega)}^p:= \sup_{\lambda > 0}\, \lambda^p \, \big| \big\{ x \in \Omega\, : \, |f(x)| \geq \lambda \big\} \big|.
\end{equation}

\subsection*{Acknowledgements.} A. J. F. and T. L. wish to thank Dr. Matteo Cozzi for very helpful discussions as well as Prof. Piotr Haj\l asz for his kind help and for pointing out some useful references. The authors also thank the anonymous referees for their valuable comments and corrections. This work has received funding from the European Research Council (ERC) under the European Union's Horizon 2020 research and innovation programme through the Consolidator Grant agreement 862342 (A. J. F.). B. A. is partially supported by projects MTM2016-80474-P and PID2019-110712GB-I00, MINECO, Spain. B. A. and A. Y. are partially supported by the DGRSDT, Algeria.

\section{Preliminary integral estimates and useful tools} \label{estimates-section}

\noindent We begin this section with an elementary (simple calculus) result. 

\begin{lemma}\label{lemma-MVT}
Let $N \geq 1$ and $\lambda > 1$, then:
$$
\big| |u|^{\lambda} - |w|^{\lambda} \big| \leq \lambda \max\{1,2^{\lambda-2} \} \big( |u-w|^{\lambda-1} + |u|^{\lambda-1}\big)  |u-w| , \quad \textup{ for all } u,w \in \RN.
$$
\end{lemma}

\begin{proof}
Let us fix arbitrary $u, w \in \RN$ and consider $f: \R \to \R$ given by $f(\tau) = |\tau u + (1-\tau) w|^{\lambda}$. Then, since $\lambda > 1$, observe that
\begin{equation*}
f(1) = |u|^{\lambda}, \quad f(0) = |w|^{\lambda}, \quad \textup{ and } \quad f \in C^1(\R).
\end{equation*}
Thus, using the Mean Value Theorem, we get
\begin{equation} \label{mvt-eq1}
\big| |u|^{\lambda} - |w|^{\lambda} \big| \leq \sup_{\tau \in [0,1]} |f'(\tau)|.
\end{equation}
On the other hand, using the Cauchy-Schwarz inequality, the triangular inequality and the fact that
$$
(a+b)^p \leq \max \{1,2^{p-1}\} (a^p + b^p), \quad \textup{ for all } a, b \geq 0 \textup{ and } p > 0,
$$
we get
\begin{equation} \label{mvt-eq2}
\begin{aligned}
|f'(\tau)| & = \lambda |\tau u + (1-\tau)w|^{\lambda-2} | \langle \tau u + (1-\tau)w, u-w \rangle |
\\
& = \lambda |u-w| |(1-\tau)(w-u) + u|^{\lambda-1}  \leq \lambda|u-w| \big( (1-\tau)|u-w| + |u| \big)^{\lambda-1} \\
& \leq \lambda \max\{1,2^{\lambda-2}\} \big( (1-\tau)^{\lambda-1} |u-w|^{\lambda-1} + |u|^{\lambda-1} \big)|u-w|,
 \qquad \textup{ for all } \tau \in [0,1].
\end{aligned}
\end{equation}
The result immediately follows combining \eqref{mvt-eq1} and \eqref{mvt-eq2}.
\end{proof}

\subsection{Preliminary integral estimates} \label{sec40} $ $ \medbreak

\noindent The aim of this subsection is to prove integral estimates which will be of key importance to the proof of our main results. The following lemma is essentially contained in the proof of \cite[Lemma 3.13]{G-R-2007}. For the benefit of the reader we give here a self-contained direct calculus proof.

\begin{lemma}\label{Grzywny} Let $N \geq 1$, $\rho > 0$ and $\alpha, \beta \in (-\infty,N)$. There exists $C:= C(N,\rho, \alpha, \beta) > 0$ such that:
\begin{itemize}
\item[$\bullet$] If $ N-\alpha - \beta\neq 0$, then
$$
\int_{B_\rho (0)} \frac{dz}{|x-z|^{\alpha} |y-z|^{\beta} }  \leq  C \Big(1+ |x-y|^{N-\alpha - \beta} \Big), \quad \textup{ for all } x,y \in B_\rho(0) \textup{ with } x \neq y.
$$
\item[$\bullet$] If $ N - \alpha - \beta = 0$, then
$$
\int_{B_\rho (0)} \frac{dz}{|x-z|^{\alpha} |y-z|^{\beta} }  \leq  C \Big(1+ \big|\ln |x-y|\big| \Big),  \quad \textup{ for all } x,y \in B_\rho(0)\textup{ with } x \neq y.
$$
\end{itemize}
\end{lemma}

\begin{remark}
Following the proof it is possible to obtain the explicit dependence of the constant $C$ with respect to $\rho$. However, since it is not necessary to our purposes and it deteriorates the presentation, we prefer not to give the explicit dependence.
\end{remark}

\begin{proof}[Proof of Lemma \ref{Grzywny}] Let us fix arbitrary $x,y \in B_\rho(0)$ with $x \neq y$. We consider the map $T: \RN \to \RN$ given by $T(z) = \frac{z-x}{|x-y|}$ and note that $T(B_\rho(0)) = B_{\frac{\rho}{|x-y|}}(\frac{-x}{|x-y|}) =: A$. Hence, it follows that
$$
\int_{B_\rho (0)} \frac{dz}{|x-z|^{\alpha} |y-z|^{\beta}}
=
|x-y|^{N-\alpha - \beta} \int_{A} \frac{dw}{|w|^{\alpha} |w+ p |^{\beta}}, \quad \textup{ with } \quad p = \frac{x-y}{|x-y|}.
$$
Since $|p|=1$, we define
$$
A_1:= A \cap B_{\frac12(0)}, \quad A_2:= A \cap \big(B_2(0) \setminus B_{\frac12(0)} \big) \quad \textup{ and } \quad A_3:= A \cap \big( \RN \setminus B_2(0) \big),
$$
and observe that
\begin{equation} \label{decomposition-Grzywny}
\int_{B_\rho (0)} \frac{dz}{|x-z|^{\alpha} |y-z|^{\beta}}
=
|x-y|^{N-\alpha - \beta} \int_{A} \frac{dw}{|w|^{\alpha} |w+ p |^{\beta}} dw = |x-y|^{N-\alpha - \beta} \sum_{i = 1}^3 \int_{A_i}  \frac{dw}{|w|^{\alpha} |w+ p |^{\beta}}.
\end{equation}
We treat the integral in each $A_i$ separately. First, observe that $w \in A_1$ implies $|w+p| \geq \frac12$. Thus, using that $\alpha < N$, we get
\begin{equation} \label{A1}
\int_{A_1} \frac{dw}{|w|^{\alpha}|w+p|^{\beta}} \leq 2^{\beta} \int_{B_\frac12(0)} \frac{dw}{|w|^{\alpha}} = \frac{2^{\beta+\alpha-N}}{N-\alpha} \sigma(\mathbb{S}^{N-1}).
\end{equation}
Next, observe that $|w| \geq \frac12$ in $A_2$ and that $B_{2}(0) \setminus B_{\frac12}(0) \subset B_4(-p)$. Thus, using that $\beta < N$, we obtain
\begin{equation*}
\int_{A_2} \frac{dw}{|w|^{\alpha}|w+p|^{\beta}} \leq 2^{\alpha} \int_{B_4(-p)} \frac{dw}{|w+p|^{\beta}} = \frac{2^{2N-2\beta+\alpha}}{N-\beta} \sigma(\mathbb{S}^{N-1}).
\end{equation*}
Finally, note that $|w+p| \geq \frac{|w|}{2}$ in $A_3$ and that $A_3 \subset \{w \in \RN : 2 \leq |w| \leq \frac{2\rho}{|x-y|} \}  $. Thus, we get
\begin{equation*}
\int_{A_3} \frac{dw}{|w|^{\alpha}|w+p|^{\beta}} \leq
2^{\beta} \sigma(\mathbb{S}^{N-1}) \int_2^{\frac{2\rho}{|x-y|}} r^{N-1-\alpha-\beta} dr.
\end{equation*}
On one hand, if $N - \alpha - \beta= 0$, we get that
\begin{equation*}
\int_{A_3} \frac{dw}{|w|^{\alpha}|w+p|^{\beta}} \leq 2^{\beta} \sigma(\mathbb{S}^{N-1}) \left( |\ln \rho | + \big| \ln|x-y| \big| \right).
\end{equation*}
On the other hand, if $N - \alpha - \beta \neq 0$, we obtain that
\begin{equation} \label{A3-non-log}
\int_{A_3} \frac{dw}{|w|^{\alpha}|w+p|^{\beta}} \leq \frac{2^{N-\alpha}}{N-\alpha-\beta}\sigma(\mathbb{S}^{N-1}) \left( 1+ \frac{\rho^{N-\alpha-\beta}}{|x-y|^{N-\alpha-\beta}}\right).
\end{equation}
The result follows gathering \eqref{A1}--\eqref{A3-non-log} and substituting in \eqref{decomposition-Grzywny}.
\end{proof}

Our next result is inspired by \cite[Lemma 2.5]{J-S-W-2020} and its proof relies on \cite[Lemma 2.3]{J-S-W-2020}. We believe this result is of independent interest and will be useful to analyze the behaviour of solutions to non-local problems in bounded domains. Note that the $C^2$--regularity imposed on the boundary of the domain $\partial \Omega$ is crucial here.

\begin{lemma} \label{lemma-Tobias}
For any  $0 < \lambda < 1$ and $0 < a < 1$, there exists $C:= C(N, \lambda, a, \Omega) > 0$ such that
$$
\int_{\Omega} |x-y|^{\lambda-N} \delta^{-a}(y)\, dy  \leq C \big(1+ \omega_{a,\lambda}(x) \big), \quad \textup{ for all } x \in \Omega,
$$
where
\begin{equation*} 
\omega_{a,\lambda} (x) :=
\begin{cases}
1, & \mbox{ if } \lambda > a, \\
|\ln \delta(x)|, \qquad &\mbox{ if }  \lambda = a, \\
\delta^{-(a-\lambda)}(x), &\mbox{ if } \lambda < a.
\end{cases}
\end{equation*}
\end{lemma}

\begin{proof}
First observe that
\begin{equation} \label{decomposition-Tobias}
\begin{aligned}
\int_{\Omega} & |x-y|^{\lambda-N} \delta^{-a}(y)\, dy \\
& = \int_{B_{\frac{\delta(x)}{2}}(x)} \frac{dy}{\delta^{a}(y) |x-y|^{N-\lambda}} + \int_{\Omega \setminus B_{\frac{\delta(x)}{2}}(x)}\frac{dy}{\delta^{a}(y) |x-y|^{N-\lambda}} =: I_1(x) + I_2(x), \qquad \textup{ for all } x \in \Omega.
\end{aligned}
\end{equation}
We estimate each $I_i$ separately. First, observe that $y \in B_{\frac{\delta(x)}{2}}(x)$ implies $\delta(y) \geq \frac{\delta(x)}{2}$. Hence, there exists $C:= C(N,\delta,a,\Omega) > 0$ such that
\begin{equation} \label{I1-Tobias}
\begin{aligned}
I_1(x) & \leq 2^a \delta^{-a}(x) \int_{B_{\frac{\delta(x)}{2}}(x)} \frac{dy}{|x-y|^{N-\lambda}}  = 2^a \sigma(\mathbb{S}^{N-1}) \delta^{-a}(x) \int_0^{\frac{\delta(x)}{2}} \rho^{\lambda-1} \\
& = \frac{2^{a-\lambda}}{\lambda} \sigma(\mathbb{S}^{N-1})\, \delta^{\lambda-a}(x) \leq C \big( 1+ \omega_{a,\lambda}(x) \big), \qquad \textup{ for all } x \in \Omega.
\end{aligned}
\end{equation}
To estimate $I_2$, note that, since $\Omega$ is of class $C^2$, there exists $\epsilon = \epsilon(\Omega) \in (0,1)$  and a unique reflection map
$$
\mathcal{R}: B_{\epsilon}(\partial \Omega) \to B_{\epsilon}(\partial \Omega), \quad x \mapsto \mathcal{R}(x) := x + 2(p(x) - x),
$$
at the boundary $\partial \Omega$. Here, $B_{\epsilon}(\partial \Omega) := \{x \in \RN: \d(x)< \epsilon\}$ and $p(x)$ is the projection of $x \in B_{\epsilon}(\partial \Omega)$ onto $\partial \Omega$, i.e. the unique point on $\partial \Omega$ with $|x-p(x)| = \d(x)$. We now consider separately two cases.
\medbreak
\textbf{Case 1 :} $x \in \Omega \setminus B_{\epsilon}(\partial \Omega)$.
\smallbreak
\noindent $\qquad$ In this case, $y \in \Omega \setminus B_{\frac{\delta(x)}{2}}(x)$ implies $|y-x| \geq \frac{\epsilon}{2}$. Hence, it follows that
\begin{equation} \label{I2-case1}
I_2(x) \leq \Big(\frac{2}{\epsilon} \Big)^{N-\lambda} \int_{\Omega} \frac{dy}{\delta^{a}(y)} \leq C,
\end{equation}
$\qquad$ with $C= C(N,\lambda,a,\Omega) > 0$.
\medbreak
\textbf{Case 2 :} $x \in \Omega \cap B_{\epsilon}(\partial\Omega)$.
\smallbreak
\noindent $\qquad$ Note that, $x \in \overline{\Omega} \cap B_{\epsilon}(\partial\Omega)$ implies that
$$\mathcal{R}(x) \in B_{\epsilon}(\partial\Omega) \setminus \Omega,\quad
\delta (x) = \delta(\mathcal{R}(x)), \quad \textup{ and } \quad |x-\mathcal{R}(x)| = 2\delta(x).
$$
$\qquad$ Hence, we have that
$$
|\mathcal{R}(x)-y| \leq 5|x-y|, \quad \textup{ for all } y \in \Omega \setminus B_{\frac{\delta(x)}{2}}(x),
$$
$\qquad$ and thus
$$
I_2(x) \leq 5^{N-\lambda} \int_{\Omega \setminus B_{\frac{\delta(x)}{2}}(x)} \frac{dy}{\delta^{a}(y)|\mathcal{R}(x)-y|^{N-\lambda}}.
$$
$\qquad$ Applying then \cite[Lemma 2.3]{J-S-W-2020} to $\mathcal{R}(x) \in   \Omega^c $, we get $C = C(N,\lambda,a,\Omega) > 0$ such that
\begin{equation} \label{I2-case2}
I_2(x) \leq C \big( 1+ \omega_{a,\lambda}(x) \big).
\end{equation}
\noindent Gathering \eqref{I2-case1}-\eqref{I2-case2}, we get the existence of $C = C(N,\lambda,a,\Omega) > 0$ such that
\begin{equation} \label{I2-Tobias}
I_2(x) \leq C \big( 1+ \omega_{\lambda,a}(x) \big), \quad \textup{ for all } x \in \Omega.
\end{equation}
The result follows substituting \eqref{I1-Tobias} and \eqref{I2-Tobias} in \eqref{decomposition-Tobias}.
\end{proof}

\subsection{Useful tools} $ $ \medbreak

\noindent In this subsection we gather some known results that will be crucially used in this paper. Let us start with a classical (but difficult to find in the literature) \textit{a.e.} fundamental Theorem of Calculus in Sobolev spaces.

\begin{lemma} \label{a.e.-fundamental-th-calculus}
Let $f \in W_{loc}^{1,1}(\RN)$,  then
$$
\textup{ for a.e. } x,z, \in \RN, \quad f(x) - f(z) = \int_0^1 \langle x-z, \nabla f(\tau x + (1-\tau)z) \rangle\, d\tau\, .
$$
\end{lemma}

\begin{proof}
The result follows from the combination of \cite[Theorem 16]{H2000} (applied with $k = p = 1$) with \cite[Theorem 4.20]{evans-garipiery} and \cite[Theorem 3.35]{folland}. See also \cite[Chapter 3, Exercise 3.15]{ziemer}.
\end{proof}

The next tool we need is a classical result from harmonic analysis. It is crucial in the proofs of our main regularity results of Section \ref{sec5}.

\begin{lemma}\label{Stein}
Let $\omega \subset \RN$, $N \geq 2$, be an open bounded domain and let $0 < \alpha < N$ and  $1\le p<\ell<\infty$ be such that  $\dfrac{1}{\ell} =\dfrac{1}{p}-\dfrac{\a}{N}$. Moreover, for $g\in L^p(\omega)$, let
$$
J_\alpha(g)(x): =\int_{\omega} \dfrac{g(y)}{|x-y|^{N-\a}}dy\,.
$$
It follows that there exists  $C= C(N,\a, p , \s > 0, \ell, \O)>0 $  such that
\begin{itemize}
\item[$a)$] $J_{\alpha}$ is well defined (in the sense that the integral converges absolutely for a.e. $x \in \omega$).
\item[$b)$] If $p = 1$, then $[J_{\alpha}(g)]_{M^{\ell}(\omega)} \leq C\, \|g\|_{L^p(\omega)}.$  In particular, $\|J_{\alpha}(g)\|_{L^{\sigma}(\omega)} \leq C \, \|g\|_{L^p(\omega)}$ for all $1 \leq \sigma < \ell$.
\item[$c)$] If $1<p<\frac{N}\a$, then $\|J_{\alpha}(g)\|_{L^\ell (\omega)} \leq C\, \|g\|_{L^p(\omega)}$.
\item[$d)$] If $p = \frac{N}{\alpha}$, then $\|J_{\alpha}(g)\|_{L^{\sigma}(\omega)} \leq C \, \|g\|_{L^p(\omega)}$ for all $1 \leq \sigma < +\infty$.
    \item[$e)$]  If $p>\frac{N}\a$, then $
\|J_{\alpha}(g)\|_{L^{\infty}(\omega)} \leq C\, \|g\|_{L^p(\omega)}.
$
\end{itemize}
\end{lemma}

\begin{proof}
Parts a), b), c) and d) follow from \cite[Theorem I, Section 1.2, Chapter V, page 119]{Stein}. Part e) is contained in \cite[Lemma 7.12]{G_T_2001_S_Ed} (see also \cite[Theorem 2.2]{mizuta}).
\end{proof}

In the next lemma we collect some known estimates on the $s$-Green kernel $G_s$ and its gradient $\nabla_x G_s$.

\begin{lemma} \label{Green-estimates-lemma}
Let  $G_s(x,y)$  be the Green Function associated to   $ \Omega$ , then there exist $\mathcal{C}_1:= \mathcal{C}_1 (N,s,\Omega) > 0$ and $\mathcal{C}_2:= N$ such that
\begin{equation}\label{green100}
G_s(x,y)\le \mathcal{C}_1\min\left\{\frac{1}{|x-y|^{N-2s}}, \frac{\delta^s(x)}{|x-y|^{N-s}}, \frac{\delta^s(y)}{|x-y|^{N-s}}\right\}, \quad \textup{ for a.e. } x, y \in \Omega,
\end{equation}
and
\begin{equation} |\nabla_x G_s(x,y)|\le \mathcal{C}_2\, G_s(x,y)\max\left\{\frac{1}{|x-y|}, \frac{1}{\delta(x)}\right\}, \quad \textup{ for a.e. } x, y \in \Omega.\label{green200}
\end{equation}
\end{lemma}

\begin{remark}
Under the assumptions of the previous lemma, sharp (two-sides) pointwise estimates for the $s$-Green function are available in the literature. See for instance \cite[Corollary 1.3]{CS1998} or \cite[Corollary 1.2]{CKS2010}.
\end{remark}

\begin{proof}
The pointwise estimates for the Green function, namely \eqref{green100}, can be found in \cite[Theorem 1.1]{CS1998}, while   \eqref{green200} follows from \cite[Corollary 3.3]{BKN2002}.
\end{proof}

The next proposition will be systematically used in this paper,  too.

\begin{proposition}\label{hardy-reg}
Let $s \in (0,1)$ and $u := \mathbb{G}_s[f]$. Then, for
\begin{equation} \label{reg-u}
\begin{cases}
r= \infty, \qquad & \mbox{ if } \ m >\frac{N}{2s}\,, \\
r\in [1, \infty),  & \mbox{ if }  \ m =\frac{N}{2s}\,, \\
r=\frac{mN}{N-2ms}, & \mbox{ if }  \ 1<m <\frac{N}{2s}\,, \\
 r\in \big[1, \frac{N}{N-2s} \big), & \mbox{ if }  \ m =1\,,
\end{cases}
\qquad \quad \mbox{ and } \qquad \quad
\begin{cases}
 q=\infty, \qquad & \mbox{ if }  \ m >\frac{N}{s}\,, \\
 q\in [1, \infty),  & \mbox{ if } \  m =\frac{N}{s}\,, \\
q= \frac{mN}{N-ms}, & \mbox{ if } \  1<m <\frac{N}{s}\,, \\
 q \in \big[1, \frac{N}{N-s}  \big), & \mbox{ if }   \ m =1\,,
\end{cases}
\end{equation}
there exists $C := C(N,s,m,r,q,\O)> 0$
such that
\begin{equation} \label{reg-u1}
\|u\|_{L^{r}(\O)}+\left\|\frac{u}{\d^s}\right\|_{L^{q}(\O)}\le C\, \|f\|_{L^m(\O)}\,.
\end{equation}
\end{proposition}

\begin{proof}
The $L^r$--regularity of the function $u$ can be found in \cite{LPPS2015}. See also \cite[Proposition 1.4]{RO-S14} for a different proof. As far as the  Hardy-type regularity is concerned, we observe that, thanks to \eqref{green100},
$$
\left|\frac{u (x)}{\d^s (x) } \right| = \frac1{\d^s (x) }  \left|\int_{\Omega} G_s(x,y) f(y) dy \right|
\leq  \mathcal{C}_1  \int_{\Omega}  \frac{| f(y)|}{|x-y|^{N-s}} dy \,,  \quad  \mbox{ for a.e. } x \in \O.
$$
Hence the $L^q$-summability of $u/\d^s$ in \eqref{reg-u}--\eqref{reg-u1} follows combining
Lemma \ref{Stein} with Lemma \ref{Green-estimates-lemma}.
\end{proof}

\section{The $\frac{t}{2}$--Laplacian and the $t$--Gradient of the $s$--Green function}\label{sec4}

\noindent This section is devoted to prove Theorems \ref{main-th-Green} and \ref{gradGreen-intro}. Moreover, we will provide pointwise estimates for the $\frac{t}{2}$--Laplacian of the $s$-Green function in $\Omega^c $.

Let us begin with an elementary technical lemma that will be key in the proof of Theorems \ref{main-th-Green} and \ref{gradGreen-intro}.

\begin{lemma} \label{lemma-alpha-y}
Let $s \in (0,1)$ and  $0 < t < \min \{1,2s\}$.  Also, for any $y \in \Omega$, let $\alpha_y : \RN \setminus \{y\} \to \R$ be given by
\begin{equation} \label{alpha-y}
\alpha_y(x) := |x-y|^{N-(2s-t)} G_s(x,y).
\end{equation}
Then
$$
\alpha_y(x) -  \alpha_y(z) =  \int_0^1 \langle x-z, \nabla \alpha_y (\tau x + (1-\tau)z) \rangle d\tau, \quad \textup{ for a.e. } x,z \in \RN.
$$
\end{lemma}

\begin{proof}
For any $y \in \Omega$, we want to prove that $\alpha_y \in W_{loc}^{1,1}(\RN)$ so that the result immediately follows from Lemma \ref{a.e.-fundamental-th-calculus}. We will actually show that $\alpha_y \in W^{1,p}(\RN)$ for all $1 \leq p < \min\big\{\frac{1}{1-s}, \frac{1}{1-t}\big\}$.
Using \eqref{green100}--\eqref{green200}, we infer that
$$
|\nabla \alpha_y(x)| \leq \frac{\C_1(N-(2s-t)+\C_2)}{|x-y|^{1-t}} + \frac{\C_1\, \C_2|x-y|^{t-s}}{\delta^{1-s}(x)}\, \mathds{1}_{\{\delta(x) \leq |x-y|\}}, \quad \textup{ for a.e. } x,y \in \Omega.
$$
We now consider separately two cases:
\medbreak
\noindent \textbf{Case 1 :} $s \leq t < \min\{1,2s\}$.
\smallbreak
In this case we have that there exists $C:= C(N,s,t,p,\Omega) > 0$ such that
$$
\int_{\Omega} |\nabla \alpha_y(x)|^p dx \leq C \left( \int_{B_{d_{\Omega}}(y)} \frac{dx}{|x-y|^{p(1-t)}} + \int_{\Omega} \frac{dx}{\delta^{p(1-s)}(x)} \right),
$$
where $d_{\Omega} := \diam (\O)$.  Thus, we immediately deduce that $\alpha_y \in W_0^{1,p}(\Omega)$, for all $1 \leq p < \frac{1}{1-s}$, and so that $\alpha_y \in W^{1,p}(\RN)$, for all $1 \leq p < \frac{1}{1-s}$.

\medbreak
\noindent \textbf{Case 2 :} $0 < t < s$.
\smallbreak
In this case we have that there exists $C:= C(N,s,t,p,\Omega) > 0$ such that
$$
\int_{\Omega} |\nabla \alpha_y(x)|^p dx \leq C \left( \int_{B_{d_{\Omega}}(y)} \frac{dx}{|x-y|^{p(1-t)}} + \int_{\Omega} \frac{dx}{\delta^{p(1-t)}(x)} \right).
$$
Thus, we immediately deduce that $\alpha_y \in W_0^{1,p}(\Omega)$, for all $1 \leq p < \frac{1}{1-t}$, and so that $\alpha_y \in W^{1,p}(\RN)$, for all $1 \leq p < \frac{1}{1-t}$.
\end{proof}

Having at hand the previous lemma we can now provide the proof of Theorem  \ref{main-th-Green}.

\begin{proof}[Proof of Theorem \ref{main-th-Green}]
First of all, observe that
\begin{equation} \label{decomposition-t-Laplacian-Green}
(-\Delta)_x^{\frac{t}{2}}G_s(x,y) = G_s(x,y) \int_{  \Omega^c } \frac{dz}{|x-z|^{N+t}} + (-\Delta)_{\Omega,x}^{\frac{t}{2}}G_s(x,y), \quad    \textup{ for a.e. } x, y \in \Omega.
\end{equation}
Now, observe that,
$$
\int_{  \Omega^c } \frac{dz}{|x-z|^{N+t}} \leq \int_{ \RN \setminus B_{\delta(x)}(x) } \frac{dz}{|x-z|^{N+t}} = \sigma(\mathbb{S}^{N-1}) \int_{\delta(x)}^{\infty} \frac{d\rho}{\rho^{1+t}} = t^{-1} \sigma(\mathbb{S}^{N-1})\, \delta^{-t}(x), \quad \textup{ for all } x \in \Omega.
$$
Hence, using that $G_s(x,y) \geq 0$ for all $(x, y) \in \R^{2N}_{\ast}$, it follows that,
$$
G_s(x,y) \int_{  \Omega^c } \frac{dz}{|x-z|^{N+t}} \leq t^{-1} \sigma(\mathbb{S}^{N-1})\, G_s(x,y)\delta^{-t}(x), \quad \textup{ for a.e. } x, y \in \Omega.
$$
In particular, taking into account \eqref{green100}, we infer that
\begin{equation} \label{I1}
\begin{aligned}
G_s(x,y) \int_{  \Omega^c } \frac{dz}{|x-z|^{N+t}} & \leq  \frac{\C_1 \, t^{-1} \sigma(\mathbb{S}^{N-1})}{|x-y|^{N-s}} \frac{1}{\delta^{t-s}(x)} , \quad \textup{ for a.e. } x, y \in \Omega.
\end{aligned}
\end{equation}
On the other hand, observe that
\begin{equation} \label{decomposition-I2}
\begin{aligned}
| x  -y &|^{N-(2s-t)}  (-\Delta)_{\Omega,x}^{\frac{t}{2}}G_s(x,y) \\
&  = (-\Delta)_{\Omega,x}^{\frac{t}{2}} \alpha_y(x) + \int_{\Omega} G_s(z,y) \frac{|z-y|^{N-(2s-t)} - |x-y|^{N-(2s-t)}}{|x-z|^{N+t}} dz, \  \quad \textup{ for a.e. } x, y \in \Omega,
\end{aligned}
\end{equation}
where we recall that $\alpha_y$ has been defined in \eqref{alpha-y}. Using that $G_s(x,y) \geq 0$ for all $(x, y) \in \R^{2N}_{\ast}$ and Lemma \ref{lemma-MVT} with $\lambda = N-(2s-t)$, $u = z-y$ and $w = x-y$, we get
$$
\begin{aligned}
\bigg|\int_{\Omega} &  G_s(z,y) \frac{|z-y|^{N-(2s-t)} - |x-y|^{N-(2s-t)}}{|x-z|^{N+t}} dz\, \bigg| \\ & \leq C  \int_{\Omega} G_s(z,y) \left( |x-z|^{-2s} +  \frac{|z-y|^{N-(2s-t)-1}}{|x-z|^{N+t-1}} \right) dz,  \quad \textup{ for a.e. } x, y \in \Omega,
\end{aligned}
$$
where $C := [N-(2s-t)] \max\{1,2^{N-(2s-t)-2}\} > 0$. Thus, taking into account \eqref{green100} and applying Lemma \ref{Grzywny} with $\rho = R$ (recall that $R$ has been chosen in \eqref{erre}), we get
\begin{equation} \label{Gs-z-final}
\begin{aligned}
\bigg|\int_{\Omega}  & G_s(z,y) \frac{|z-y|^{N-(2s-t)} - |x-y|^{N-(2s-t)}}{|x-z|^{N+t}} dz \, \bigg|\\
& \leq C \left( \int_{\Omega} \frac{dz}{|x-z|^{2s} |y-z|^{N-2s}} + \int_{\Omega} \frac{dz}{|x-z|^{N+t-1}|y-z|^{1-t}} \right) \\
& \leq C \left( \int_{B_R(0)} \frac{dz}{|x-z|^{2s} |y-z|^{N-2s}} + \int_{B_R(0)} \frac{dz}{|x-z|^{N+t-1}|y-z|^{1-t}} \right) \\
& \leq C \left( 1+ \left| \ln \frac{1}{|x-y|} \right| \right), \quad \textup{ for a.e. } x, y \in \Omega.
\end{aligned}
\end{equation}
Here, and in the rest of the proof, we denote by
$C$ any constant depending only on $N,s,t $ and $\Omega$ and whose value may change from line to line.

Next, we deal with $(-\Delta)_{\Omega,x}^{\frac{t}{2}}\alpha_y(x)$. First, using Lemma \ref{lemma-alpha-y}, we write
$$
\begin{aligned}
(-\Delta)_{\Omega,x}^{\frac{t}{2}} \alpha_y (x) & = \int_{\Omega} \left( \int_0^1 \langle x-z, \nabla \alpha_y (\tau x + (1-\tau) z) \rangle d\tau \right) \frac{dz}{|x-z|^{N+t}}, \quad \textup{ for a.e. } x, y \in \Omega.
\end{aligned}
$$
Then, let us introduce the shortened notation $\theta_{\tau} := \tau x + (1-\tau)z$ and define $I_{\Omega}:= \{\tau \in [0,1]: \theta_{\tau} \in \Omega\} \subset [0,1]$. By direct computations, taking into account \eqref{green100}--\eqref{green200}, we get
\begin{equation} \label{decomposition-reg-Lap-alpha-y}
\begin{aligned}
|(-\Delta)_{\Omega,x}^{\frac{t}{2}} \alpha_y(x)| & \leq \int_{\Omega} \left( [N-(2s-t)] \int_0^1 |\theta_\tau-y|^{N-(2s-t)-1} G_s(\theta_{\tau},y) d\tau \right. \\
& \qquad\ \left. + \int_0^1 |\theta_\tau-y|^{N-(2s-t)} |\nabla_x G_s(\theta_{\tau},y)| d\tau \right) \frac{dz}{|x-z|^{N+t-1}} \\
& \leq  \C_1 \big( [N-(2s-t)] + \C_2\big)\int_{\Omega}\left( \int_{I_{\Omega}} |\theta_\tau-y|^{-(1-t)} d\tau \right) \frac{dz}{|x-z|^{N+t-1}} \\
& \qquad\ +\C_1\, \C_2\int_{\Omega} \left(  \int_{I_{\Omega}} \frac{|\theta_{\tau}-y|^{t-s}}{\delta^{1-s} (\theta_{\tau})} d\tau \right) \frac{dz}{|x-z|^{N+t-1}}\\
& =: \C_1 \big( [N-(2s-t)] + \C_2\big)  J_1(x,y) +  \C_1 \, \C_2\, J_2(x,y),  \quad \textup{ for a.e. } x, y \in \Omega.
\end{aligned}
\end{equation}
We treat each $J_i$ for separately. Thanks to the choice of $R$ (see \eqref{erre}) we have that $\Omega \subset B_R(0)$ and that $\theta_{\tau} \in B_R(0)$ for all $\tau \in [0,1]$. Thus, using Lemma \ref{Grzywny}, we get that
\begin{equation} \label{J1-final}
\begin{aligned}
J_1(x,y) & \leq \int_{B_{R}(0)} \left( \int_{0}^1 |\theta_\tau-y|^{-(1-t)} d\tau \right) \frac{dz}{|x-z|^{N+t-1}} = \int_0^1 \left( \int_{B_R(0)} \frac{dz}{|y-\theta_{\tau}|^{1-t}|x-z|^{N+t-1}} \right) d\tau  \\
& = \int_0^1 \left( \int_{B_R(0)} \frac{dz}{|(1-\tau)^{-1}(y-\tau x) - z|^{1-t} |x-z|^{N+t-1}} \right) \frac{d\tau}{(1-\tau)^{1-t}}  \\
& \leq C \int_0^1 \frac{1}{(1-\tau)^{1-t}} \left( 1+ \left| \ln \left( \frac{|x-y|}{1-\tau} \right) \right| \right) d\tau  \leq C \big(1+\big|\ln {|x-y|} \big|\big), \quad \textup{ for a.e. } x, y \in \Omega.
\end{aligned}
\end{equation}
We now deal with $J_2$. First, let us define $\Omega_{\tau,x}:= \{z \in \Omega: \tau x + (1-\tau)z \in \Omega\}$ and $T(z) = \theta_{\tau} = \tau x + (1-\tau)z$, and observe that $T(\Omega_{\tau,x}) \subseteq \Omega$. Hence, using the change of variable $\theta = T(z)$, we get
$$
\begin{aligned}
J_2(x,y) & = \int_0^1 \left( \int_{\Omega_{\tau,x}} \frac{|\theta_{\tau}-y|^{t-s}}{\delta^{1-s}(\theta_{\tau})} \frac{dz}{|x-z|^{N+t-1}} \right) d\tau = \int_0^1 \left( \int_{T(\Omega_{\tau,x})} \frac{|\theta-y|^{t-s}}{\delta^{1-s}(\theta)} \frac{(1-\tau)^{t-1}}{|\theta-x|^{N+t-1}} d \theta \right) d\tau \\
& \leq \int_0^1 (1-\tau)^{t-1} d\tau \left( \int_{\Omega} \frac{|\theta-y|^{t-s}}{\delta^{1-s}(\theta)} \frac{ d \theta }{|\theta-x|^{N+t-1}}  \right)  = \frac{1}{t} \int_{\Omega}  \frac{|\theta-y|^{t-s}}{\delta^{1-s}(\theta)} \frac{d\theta}{|\theta-x|^{N+t-1}}, \ \, \textup{for a.e. } x,y \in \Omega.
\end{aligned}
$$
Moreover, notice that $|\theta -y|^{t-s} \leq |\theta-x|^{t-s} + |x-y|^{t-s}$ for all $x, y, \theta \in \RN$. Thus, we actually have that
$$
\begin{aligned}
J_2(x,y)    \leq   \frac{1}{t} \int_{\Omega}  \frac{1}{\delta^{1-s}(\theta)} \frac{d\theta}{|\theta-x|^{N+s-1}}
 + \frac{|x-y|^{t-s}}{t}   \int_{\Omega}  \frac{d\theta}{|\theta-x|^{N+t-1}\delta^{1-s}(\theta)},
 \quad \textup{ for a.e. } x,y \in \Omega.
\end{aligned}
$$
Then, applying Lemma \ref{lemma-Tobias} with $\lambda = 1-s$ and $a = 1-s$ in the first integral and with $\lambda = 1-t$ and $a = 1-s$ in the second one, we get
\begin{equation} \label{J2-final}
J_2(x,y) \leq C \bigg( 1+\frac{|x-y|^{t-s}}{\d^{t-s} (x)} + |\log  \d (x)| \bigg), \quad \textup{ for a.e. } x,y \in \Omega.
\end{equation}
Having at hand \eqref{J1-final} and \eqref{J2-final}, we plug them into \eqref{decomposition-reg-Lap-alpha-y} and infer that
\begin{equation} \label{reg-Lap-alpha-y-final}
\big|(-\Delta)_{\Omega,x}^{\frac{t}{2}} \alpha_y (x) \big| \leq C  \left(  \big| \log|x-y| \big| +  |\log  \d (x)| +\frac{|x-y|^{t-s}}{\d^{t-s} (x)}   \right), \quad \textup{ for a.e. } x, y \in \Omega.
\end{equation}
Consequently, gathering \eqref{decomposition-I2}, \eqref{Gs-z-final} and \eqref{reg-Lap-alpha-y-final}, we obtain that
\begin{equation*} \label{reg-Lap-Green-final}
\big|(-\Delta)^{\frac{t}{2}}_{\Omega,x} \ G_s(x,y)\big|  \leq \frac{C}{|x-y|^{N-(2s-t)}} \left(  \big| \log|x-y| \big| +  |\log  \d (x)| +\frac{|x-y|^{t-s}}{\d^{t-s} (x)}   \right), \quad \textup{for a.e.\, }x, y \in \Omega.
\end{equation*}
The result follows combining the above inequality with \eqref{decomposition-t-Laplacian-Green} and \eqref{I1}.
\end{proof}

\begin{remark}
Modifying the proof of Theorem \ref{main-th-Green}, we can cover also the case $0 < t < s$. Indeed, to cover this case, we just have to modify \eqref{J2-final}. Note that, in the case where $0 < t < s$,
$$
J_2 (x,y) \leq   \frac{1}{t} \int_{\Omega}  \frac{1}{ |\theta-y|^{s-t} \ \delta^{1-s}(\theta)} \frac{d\theta}{|\theta-x|^{N+t-1}}, \quad \textup{ for a.e. } x,y \in \Omega.
$$
Then, splitting $\O$ into $\{ \theta \in \O \  : \ |\theta-y| \geq |\theta-x|\}=: \O_1$ and $\{ \theta \in \O \  : \ |\theta-y| \geq |\theta-x|\}=: \O_2$, we obtain
$$
\begin{aligned}
J_2 (x,y) & \leq   \frac{1}{t} \int_{\Omega_1}  \frac{1}{ \delta^{1-s}(\theta)} \frac{d\theta}{|\theta-x|^{N+s-1}}
+
 \frac{1}{t} \int_{\Omega_2}  \frac{1}{ \delta^{1-s}(\theta)} \frac{d\theta}{|\theta-y|^{N+s-1}}  \\
  & \leq   \frac{1}{t} \int_{\Omega}  \frac{1}{ \delta^{1-s}(\theta)} \frac{d\theta}{|\theta-x|^{N+s-1}}
+
 \frac{1}{t} \int_{\Omega}  \frac{1}{ \delta^{1-s}(\theta)} \frac{d\theta}{|\theta-y|^{N+s-1}},  \quad \textup{ for a.e. }x,y \in \Omega.
 \end{aligned}
$$
Applying twice Lemma \ref{lemma-Tobias} with $\lambda = 1-s$ and $a = 1-s$, we get a constant  $C:= C(N,s,t,\Omega) > 0$ such that
\begin{equation*}
J_2(x,y) \leq C \, \bigg( 1+ |\log  \d (x)|  + |\log  \d (y)| \bigg), \quad \textup{ for a.e. } x,y \in \Omega.
\end{equation*}
Then, arguing as in the proof of Theorem \ref{main-th-Green}, one can easily conclude that, if $0 < t < s$, there exists a constant $C:= C(N,s,t,\Omega) > 0$ such that, for almost every $x, y \in \Omega$,
\begin{equation*}
\big|(-\Delta)_x^{\frac{t}{2}} G_s(x,y)\big| \leq \frac{C}{|x-y|^{N-(2s-t)}} \left( \big|\log |x-y|\big| + |\log\delta(x)| + |\log \d(y)| + \frac{\d^{s-t} (x)}{ |x-y|^{s-t}} \right).
\end{equation*}
We believe that  the estimate we have  obtained in the case $0<t<s$ is far from being sharp when $x$ approaches the boundary of $\O$.
\end{remark}

As already announced, we also obtain pointwise estimates for the $\frac{t}{2}$-Laplacian of $G_s$ in $\Omega^c $.

\begin{proposition}  \label{Green-Out}
Let $s \in (0,1)$, $s \leq t < \min \{1,2s\}$ and $R$ as in \eqref{erre}.
There exists $C:= C(N,s,t,\Omega) >  0$ such that
\begin{equation} \label{Green-Out-close}
\big|(-\Delta)_x^{\frac{t}{2}} G_s(x,y) \big| \leq \frac{C}{|x-y|^{N-2s}} \left( \frac{1}{|x-y|^{t}} + \frac{1}{\d^t(x)} \right), \quad \textup{ for all } y \in \Omega,\ x \in B_{R}(0) \setminus \Omega.
\end{equation}
and
\begin{equation} \label{Green-Out-far}
\big|(-\Delta)_x^{\frac{t}{2}} G_s(x,y) \big| \leq \frac{C}{(1+|x|)^{N+t}}, \quad \textup{ for all } y \in \Omega,\ x \in B_{R}^c(0).
\end{equation}
\end{proposition}

\begin{proof}
First of all, observe that
$$
\big|(-\Delta)_x^{\frac{t}{2}} G_s(x,y) \big| = \int_{\Omega} \frac{G_s(z,y)}{|x-z|^{N+t}} dz, \quad \textup{ for all } y \in \Omega,\ x \in   \Omega^c .
$$
Now, according to the choice of $R$, we notice that $|x-z| \geq \frac{1}{4}(1+|x|)$ for all $x \in \RN \setminus B_{R}(0)$ and all $z \in \Omega$. Hence, taking into account \eqref{green100}, we get the existence of $C= C(N,s,t,\Omega)> 0$ such that
\begin{equation} \label{dist-geq-1}
\big|(-\Delta)_x^{\frac{t}{2}} G_s(x,y) \big| \leq \frac{2^{2(N+t)}}{(1+|x|)^{N+t}} \int_{\Omega} G_s(z,y) dz \leq \frac{C}{(1+|x|)^{N+t}}, \quad \textup{ for all } y \in \Omega,\ x \in  B_{R}^c(0).
\end{equation}
On the other hand, for $x\in B_{R}(0) \setminus \Omega$,  we split
$\Omega$ into
$$
\Omega_1 = \Big\{z \in \Omega: |y-z| > \frac12 |x-y| \Big\}\qquad \mbox{ and} \qquad
\Omega_2 = \Big\{z \in \Omega: |y-z| \leq \frac12 |x-y| \Big\},
$$
and we note that
 $$
 |x-z| \geq  |y-x|-  |y-z| \geq \frac12 |x-y|,
 \quad \mbox{ for all  } z  \in \Omega_2.
 $$
Thus, using \eqref{green100}, we get
\begin{equation*}
\begin{aligned}
\big|(-\Delta)_x^{\frac{t}{2}}& G_s(x,y) \big|
\leq \C_1
 \int_{\Omega}  \frac{1}{|z-x|^{N+t}}
 \frac{dz}{|z-y|^{N-2s}}  \\
& \leq
\frac{  2^{N-2s}\, \C_1}{|x-y|^{N-2s}}
  \int_{\Omega_1}  \frac{dz}{|z-x|^{N+t}}
 +
 \frac{  2^{N+t}\, \C_1}{|x-y|^{N+t}} \int_{\Omega_2}  \frac{dz}{|z-y|^{N-2s}}, \quad \textup{ for all } y \in \Omega,\, x \in B_{R}(0) \setminus \Omega.
\end{aligned}
\end{equation*}
Now, observe that
$$
\begin{aligned}
 \int_{\Omega_1}  \frac{dz}{|z-x|^{N+t}}
& \leq \sigma(\mathbb{S}^{N-1})
 \int_{\d(x)}^{R}  \frac{d\rho }{\rho^{1+t}}  \\
 & = \frac1{t}\sigma(\mathbb{S}^{N-1}) \bigg[  \d^{-t} (x)   - R^{-t}\bigg] \leq \frac{C }{  \d^t (x)}, \quad \textup{ for all } x \in B_{R}(0) \setminus \Omega,
\end{aligned}
$$
and
$$
\int_{\Omega_2}   \frac{dz}{|z-y|^{N-2s}}   \leq  \sigma(\mathbb{S}^{N-1})
\int_0^{\frac12 |x-y|}
  \frac1{\rho^{1-2s}} = \sigma(\mathbb{S}^{N-1}) \frac{1}{2^{1+2s}s} |x-y|^{2s}, \quad \textup{ for all } y \in \Omega,\, x \in B_{R}(0) \setminus \Omega.
$$
Hence, there exists $C:= C(N,s,t,\Omega) > 0$ such that
\begin{equation}\label{dist-leq-1}
\big|(-\Delta)_x^{\frac{t}{2}} G_s(x,y) \big| \leq \frac{C}{|x-y|^{N-2s}} \left( \frac{1}{|x-y|^{t}} + \frac{1}{\d^t(x)} \right), \quad \textup{ for all } y \in \Omega,\ x \in B_{R}(0) \setminus \Omega.
\end{equation}
The result follows from \eqref{dist-geq-1} and \eqref{dist-leq-1}.
\end{proof}

We now prove Theorem \ref{gradGreen-intro}. As mentioned in the introduction, this result follows by modifying in a suitable way the proof of Theorem \ref{main-th-Green}. For the benefit of the reader we provide some details.

\begin{proof}[Proof of Theorem \ref{gradGreen-intro}]
First of all, we notice that
\begin{equation} \label{dec-t-grad-Green}
\n_x^t G_s(x,y)
= G_s(x,y) \int_{  \Omega^c } \frac{x-z}{|x-z|^{N+t+1}}  dz
+
\int_{  \Omega  } \frac{(x-z)(G_s(x,y) - G_s(z,y))}{|x-z|^{N+t+1}} dz =: \I_1 + \I_2,
\end{equation}
for a.e. $x, y \in \Omega$. The first integral can be treated exactly as before, so that we deduce that
\begin{equation*}
\begin{aligned}
|\I_1 (x,y)| \leq G_s(x,y) \int_{  \Omega^c } \frac{dz}{|x-z|^{N+t}} & \leq  \frac{\C_1 \, t^{-1} \sigma(\mathbb{S}^{N-1})}{|x-y|^{N-s}}   \frac{1}{\delta^{t-s}(x)}, \quad \textup{ for a.e. } x, y \in \Omega.
\end{aligned}
\end{equation*}
As far as $\I_2$ is concerned, we observe that, for a.e. $x, y, \in  \Omega$, it follows that
\begin{equation} \label{dec-I2}
\begin{aligned}
| x  -y & |^{N-(2s-t)}  \I_2 (x,y) \\
&  =
\int_{  \Omega  } \frac{(x-z)(\a_y (x)- \a_y (z)) }{|x-z|^{N+t+1}} dz
+ \int_{\Omega} (x-z)G_s(z,y) \frac{|z-y|^{N-(2s-t)} - |x-y|^{N-(2s-t)}}{|x-z|^{N+t+1}}  dz,
\end{aligned}
\end{equation}
where we recall that $\alpha_y (\cdot) $ has been defined in \eqref{alpha-y}. On one hand, arguing as in the proof of \eqref{Gs-z-final}, we get $C > 0$ such that
\begin{equation*}
\begin{aligned}
\bigg|\int_{\Omega} & (x-z) G_s(z,y) \frac{|z-y|^{N-(2s-t)} - |x-y|^{N-(2s-t)}}{|x-z|^{N+t+1}} dz\, \bigg| \leq C \left( 1+ \left| \ln \frac{1}{|x-y|} \right| \right), \quad \textup{ for a.e. } x, y \in \Omega.
\end{aligned}
\end{equation*}
Here, and in the rest of the proof, we denote by $C$ any constant depending only on $N,s,t $ and $\Omega$ and whose value may change from line to line. On the other hand, to deal with the first integral in
\eqref{dec-I2}, we apply  Lemma \ref{lemma-alpha-y} and obtain that, for a.e. $x, y, \in \Omega$,
$$
\begin{aligned}
\int_{  \Omega  } \frac{(x-z)(\a_y (x)- \a_y (z)) }{|x-z|^{N+t+1}} dz
& = \int_{\Omega} (x-z) \left( \int_0^1 \langle x-z, \nabla \alpha_y (\tau x + (1-\tau) z) \rangle d\tau \right) \frac{dz}{|x-z|^{N+t+1}}  
\end{aligned}
$$
With the same notation used in the proof of Theorem \ref{main-th-Green} and thanks to   \eqref{green100}--\eqref{green200}, we get that
\begin{equation*}
\begin{aligned}
\bigg|\int_{  \Omega  } \frac{(x-z)(\a_y (x)- \a_y (z)) }{|x-z|^{N+t+1}} dz\bigg| & \leq \int_{\Omega} \left( [N-(2s-t)] \int_0^1 |\theta_\tau-y|^{N-(2s-t)-1} G_s(\theta_{\tau},y) d\tau \right. \\
& \qquad\ \left. + \int_0^1 |\theta_\tau-y|^{N-(2s-t)} |\nabla_x G_s(\theta_{\tau},y)| d\tau \right) \frac{dz}{|x-z|^{N+t-1}} \\
& \leq  \C_1 \big( [N-(2s-t)] + \C_2\big)\int_{\Omega}\left( \int_{I_{\Omega}} |\theta_\tau-y|^{-(1-t)} d\tau \right) \frac{dz}{|x-z|^{N+t-1}} \\
& \qquad\ +\C_1\, \C_2\int_{\Omega} \left(  \int_{I_{\Omega}} \frac{|\theta_{\tau}-y|^{t-s}}{\delta^{1-s} (\theta_{\tau})} d\tau \right) \frac{dz}{|x-z|^{N+t-1}}\\
& =: \C_1 \big( [N-(2s-t)] + \C_2\big)  J_1(x,y) +  \C_1 \, \C_2\, J_2(x,y),  \quad \textup{ for a.e. } x, y \in \Omega.
\end{aligned}
\end{equation*}
Using the estimates \eqref{J1-final} and \eqref{J2-final} for $J_1$ and $J_2$ respectively,  we conclude the proof as in Theorem \ref{main-th-Green}.
\end{proof}

\begin{remark}\label{f1.6}
Note that \eqref{grest} is crucial to write formula \eqref{representation-formula-t-half-laplacian}. Indeed, it guarantees that, when writing explicitly the fractional Laplacian of $u$, we can charge the fractional derivatives to the Green function.  
\end{remark}

\section{The fractional Poisson problem: global fractional Calder\'on--Zygmund type regularity}\label{sec5}

\noindent In this section we apply the pointwise estimates for the $\frac{t}{2}$-Laplacian of the $s$-Green kernel $G_s$ (namely, Theorem \ref{main-th-Green}) to obtain (global) Calder\'on--Zygmund type regularity results for the fractional Poisson problem \eqref{main0}. As explained in the introduction, the main idea in the proofs of these results relies on the use of a suitable representation formula. Let us begin with a preliminary technical (but very useful) lemma.


\begin{lemma} \label{declem}
Let $s \in (0,1)$, $s \leq t < \min\{1,2s\}$ and $u:= \mathbb{G}_s[f]$ with $f \in L^m(\Omega)$ for some $m \geq m_{\star}(s,t)$. Then, there exists $C := C(N, s, t,  \Omega) > 0$  such that
\begin{equation*} 
\big| (-\Delta)^{\frac{t}{2}} u (x)| \leq C \left[ g_1(x) + |\log \delta(x)| g_2(x) + \frac{1}{\delta^{t-s}(x)} g_3(x) \right], \quad \mbox{ for a.e. } x \in \Omega.
\end{equation*}
Here:
\begin{align} \label{g1}
& \bullet\ {g}_1 (x) := \io \frac{\big| \log |x-y|\,\big|}{|x-y|^{N-(2s-t)}} |f(y)| \  dy, \\
& \bullet\ g_2(x) := J_{2s-t}(|f|)(x) = \int_{\Omega} \frac{|f(y)|}{|x-y|^{N-(2s-t)}} dy \label{g2} , \\
& \bullet\  g_3  (x) :=  J_{s}(|f|)(x) =  \int_{\Omega} \frac{|f(y)|}{|x-y|^{N-s}} dy.  \label{g3}
\end{align}
\end{lemma}

\begin{proof}
The result is an immediate consequence of the representation formula \eqref{representation-formula-t-half-laplacian} and Theorem \ref{main-th-Green}.
\end{proof}


\begin{remark}
Let us recall some classical results about the distance function that will be very useful hereafter. Let $\Omega \subset \RN$, $N \geq 2$, be a bounded domain with boundary $\partial\Omega$ of class $C^2$. It is well-known that:
\begin{align} \label{dist}
& \bullet\ \textup{For all } \mu \in (0,1),\ \d^{-\mu} \in L^{\rho}(\Omega) \textup{ if and only if } 1 \leq \rho < \frac{1}{\mu}\ ;  \\
&\bullet\  \log \delta \in L^{\rho}(\Omega) \textup{ if and only if } 1 \leq \rho < +\infty. \label{dist2}
\end{align}
\end{remark}

The first result that we prove turns out to be  a generalization of Theorem \ref{calderon-zygmund-very-integrable-intro}.

\begin{theorem} \label{calderon-zygmund-very-integrable}
Let $s,t \in (0,1)$ and $u:= \mathbb{G}_s[f]$ with $f \in L^{m}(\Omega)$ for some $m > \max\big\{\frac{N}{s},\frac{N}{2s-t}\big\}$. Then: \smallbreak
\begin{itemize}
\item[i)] If $0 < t < s$, for all $1 \leq p \leq +\infty$, there exists $C > 0$ such that
$$
\|(-\Delta)^{\frac{t}{2}}u\|_{L^p(\RN)} \leq C \|f\|_{L^{m}(\Omega)}.
$$
\item[ii)] If $ t = s$, for all $1 \leq p < + \infty$, there exists $C > 0$ such that
$$
\|(-\Delta)^{\frac{s}{2}}u\|_{L^p(\RN)} \leq C \|f\|_{L^{m}(\Omega)}.
$$
\item[iii)] If $s < t < \min\{1,2s\}$,  for all $1 \leq p < \frac{1}{t-s}$, there exists $C > 0$ such that
$$
\|(-\Delta)^{\frac{t}{2}}u\|_{L^p(\RN)} \leq C \|f\|_{L^{m}(\Omega)}.
$$
\end{itemize}
Here, $C > 0$ are constants depending only on $N$, $s$, $t$, $p$, $m$ and $\Omega$.
\end{theorem}

\begin{proof}[Proof of Theorem \ref{calderon-zygmund-very-integrable}]
Let us start with two preliminary observations that will be very useful in the rest of the proof. First, since $m > \max\big\{\frac{N}{s}, \frac{N}{2s-t} \big\}$, by \cite[Proposition 1.4, (iii)]{RO-S14}, we know that there exists a positive constant $C_1= C_1 (N,s,m, \Omega)  > 0$  such that
\begin{equation} \label{u-Cs}
\|u\|_{C^s(\RN)} \leq C_1 \|f\|_{L^m(\Omega)}\,.
\end{equation}
Moreover, for any $0<t<\min\{1,2s\}$, we can easily estimate the $\frac{t}{2}$-Laplacian of $u$ \lq\lq far away\rq\rq{ } from $\Omega$. Indeed, thanks to the choice of $R$ (cf. \eqref{erre} and \eqref{erre2}), we have that $|x-y| \geq \frac14(1+|x|)$ for all $x \in \RN \setminus B_R(0)$ and all $y \in \Omega$. Thus, using   \eqref{u-Cs} and the definition of $(-\Delta)^{\frac{t   }{2}}$, we get $C := C(N,s,t,m,\Omega) > 0$ such that
\begin{equation*} 
\big| (-\Delta)^{\frac{t}{2}} u(x) \big| = \int_{\Omega} \frac{|u(y)|}{|x-y|^{N+t}} dy \leq  C \,   \frac{\|f\|_{L^m(\Omega)}}{(1+|x|)^{N+t}}, \quad \textup{ for a.e. } x \in \RN \setminus B_R(0).
\end{equation*}
Hence, for all $1 \leq p \leq \infty$, there exists $C:= C(N,s,t,m,p,\Omega)> 0$ such that
\begin{equation} \label{far2-2}
\| (-\Delta)^{\frac{t}{2}} u\|_{L^{p}(B_R^c(0))} \leq C \|f\|_{L^m(\Omega)}.
\end{equation}

We now consider separately three cases.

\medbreak
\noindent \textbf{Case 1 :} $0 < t < s$.
\smallbreak
Taking into account \eqref{u-Cs}, we get $C = C(N,s,t,m,\Omega)> 0$ such that
\begin{align*}
\big| (-\Delta)^{\frac{t}{2}}u(x) \big| & \leq \int_{\RN} \frac{|u(x)-u(y)|}{|x-y|^{N+t}} = \int_{B_1(x)}  \frac{|u(x)-u(y)|}{|x-y|^{N+t}} dy + \int_{\RN \setminus B_1(x)}  \frac{|u(x)-u(y)|}{|x-y|^{N+t}} dy \\
& \leq C_1 \|f\|_{L^m(\Omega)} \int_{B_1(x)} \frac{dy}{|x-y|^{N-(s-t)}} + 2 \, C_1 \|f\|_{L^m(\Omega)} \int_{\RN \setminus B_1(x)} \frac{dy}{|x-y|^{N+t}} \\
& = C_1 \sigma(\mathbb{S}^{N-1}) \|f\|_{L^m(\Omega)} \bigg( \int_0^1 \frac{d\rho}{\rho^{1-(s-t)}} + 2 \int_1^{\infty} \frac{d\rho}{\rho^{1+t}} \bigg) = C \|f\|_{L^m(\Omega)}, \quad \textup{ for a.e. } x \in B_R(0).
\end{align*}
Combining the above estimate with \eqref{far2-2}, the proof of $i)$ is concluded.
\medbreak
\noindent \textbf{Case 2 :} $t = s$.
\smallbreak
In order to obtain the estimate in $\O$, we apply  Lemma \ref{declem} and obtain that
\begin{equation*}
\big| (-\Delta)^{\frac{t}{2}} u (x)| \leq C \Big[ g_1(x) + \big(1+ |\log \delta(x)|\big) g_2(x) \Big], \quad \mbox{ for a.e. } x \in \Omega.
\end{equation*}
Since $m > \frac{N}{s}$, by Lemma \ref{Stein}, e), we get $C := C(N,s,m,\O)  > 0$
$$
g_1(x) + (1+|\log \delta(x)|) g_2(x) \leq C \|f\|_{L^m(\Omega)} \big(1+ |\log \delta(x)|\big), \quad \textup{ for all } x \in \Omega.
$$
Thus, recalling \eqref{dist2}, we infer that, for all $1 \leq p < \infty$, there exists $C:= C(N,s,m,p,\O) >0$
such that
\begin{equation} \label{frac-s-2-high-regularity-interior}
\|(-\Delta)^{\frac{s}{2}}u\|_{L^p(\Omega)} \leq C \|f\|_{L^m(\Omega)}.
\end{equation}
Finally, using Lemma \ref{hardy-reg} and the fact that $|x-y| \geq \max\{ \delta(y), \delta(x)\}$ for all $y \in \Omega$ and $x \in B_R(0) \setminus \Omega$, we obtain $C_1,C_2,C_3 > 0$ (depending only on $N$, $s$, $m$ and $\Omega$) such that
\begin{equation} \label{final---}
\begin{aligned}
\big| (-\Delta)^{\frac{s}{2}} u(x) \big| & \leq \int_{\Omega} \frac{|u(y)|}{|x-y|^{N+s}} \leq \int_{\Omega} \frac{|u(y)|}{\delta^s(y)} \frac{dy}{|x-y|^{N}} \leq C_1 \|f\|_{L^m(\Omega)} \int_{\Omega} \frac{dy}{|x-y|^{N}} \\
& \leq C_2 \|f\|_{L^m(\Omega)} \int_{\d(x)}^{2R} \frac{d\rho}{\rho} \leq C_3 \|f\|_{L^m(\Omega)} \big( 1+ |\log\delta(x)| \big), \quad \textup{ for a.e. } x \in B_R(0) \setminus \Omega.
\end{aligned}
\end{equation}
It is then clear that, for all $1 \leq p < \infty$, there exists $C = C(N,s,t,m,\Omega)> 0$ such that
\begin{equation} \label{frac-s-2-high-regularity-exterior-close}
\|(-\Delta)^{\frac{s}{2}}u\|_{L^p(B_R(0) \setminus \Omega)} \leq C \|f\|_{L^m(\Omega)}.
\end{equation}
Gathering \eqref{far2-2}, \eqref{frac-s-2-high-regularity-interior} and \eqref{frac-s-2-high-regularity-exterior-close} we get $ii)$.
\medbreak
\noindent
\textbf{Case 3 :} $s < t < \min\{1,2s\}$
\smallbreak
By Lemma \ref{declem} we know that
$$
\big| (-\Delta)^{\frac{t}{2}} u (x)| \leq C \left[ g_1(x) + |\log \delta(x)| g_2(x) + \frac{1}{\delta^{t-s}(x)} g_3(x) \right], \quad \mbox{ for a.e. } x \in \Omega.
$$
On the other hand, since $m > \frac{N}{2s-t}$, once again by Lemma \ref{Stein}, e), we know there exists $ C := C(N,s,t,m,\Omega)> 0> 0$ such that
$$
g_1(x) + |\log\d(x)| g_2(x) + \frac{1}{\d^{t-s}(x)} g_3(x)  \leq C \|f\|_{L^m(\Omega)} \Big( 1+ |\log\d(x)| + \frac{1}{\d^{t-s}(x)} \Big), \quad \textup{ for all } x \in \Omega.
$$
Thus, recalling \eqref{dist}, we infer that, for all $1 \leq p < \frac{1}{t-s}$, there exists $C := C(N,s,t,m,p,\Omega)> 0$  such that
\begin{equation} \label{frac-t-2-high-regularity-interior}
\|(-\Delta)^{\frac{t}{2}}u\|_{L^p(\Omega)} \leq C \|f\|_{L^m(\Omega)}.
\end{equation}
Finally, using Lemma \ref{hardy-reg} and the fact that $|x-y| \geq \max\{ \delta(y), \delta(x)\}$ for all $y \in \Omega$ and $x \in B_R(0) \setminus \Omega$, we obtain $C_1,C_2,C_3 > 0$ (depending only on $N$, $s$, $t$, $m$ and $\Omega$) such that
\begin{equation} \label{out-weight-high-integrability}
\begin{aligned}
\big| (-\Delta)^{\frac{t}{2}} u(x) \big| & \leq \int_{\Omega} \frac{|u(y)|}{|x-y|^{N+t}} \leq \int_{\Omega} \frac{|u(y)|}{\delta^s(y)} \frac{dy}{|x-y|^{N+(t-s)}} \leq C_1 \|f\|_{L^m(\Omega)} \int_{\Omega} \frac{dy}{|x-y|^{N+(t-s)}} \\
& \leq C_2 \|f\|_{L^m(\Omega)} \int_{\d(x)}^{2R} \frac{d\rho}{\rho^{1+(t-s)}} \leq C_3 \|f\|_{L^m(\Omega)} \big( 1+ \frac{1}{\d^{t-s}(x)} \big), \quad \textup{ for a.e. } x \in B_R(0) \setminus \Omega.
\end{aligned}
\end{equation}
It then follows that, for all $1 \leq p < \frac{1}{t-s}$, there exists $C = C(N,s,t,m,p,\Omega)> 0$   such that
\begin{equation} \label{frac-t-2-high-regularity-exterior-close}
\|(-\Delta)^{\frac{t}{2}}u\|_{L^p(B_R(0) \setminus \Omega)} \leq C \|f\|_{L^m(\Omega)}.
\end{equation}
Gathering \eqref{far2-2},  \eqref{frac-t-2-high-regularity-interior} and \eqref{frac-t-2-high-regularity-exterior-close} the result follows also in this case.
\end{proof}

\begin{corollary} \label{cor1-high-integrability}
Let $s, t \in (0,1)$ and $u:= \mathbb{G}_s[f]$ with $f \in L^{m}(\Omega)$ for some $m > \max\big\{\frac{N}{s},\frac{N}{2s-t}\big\}$:
\begin{itemize}
\item[i)] If $0 < t \leq s$, then $u \in L^{t,p}(\RN)$ for all $1 < p < \infty$ and
$
\|u\|_{L^{t,p}(\RN)} \leq C \|f\|_{L^m(\Omega)}.
$
\item[ii)] If $s < t < \min\{1,2s\}$, then $u \in L^{t,p}(\RN)$ for all $1 <p< \frac{1}{t-s}$ and
$
\|u\|_{L^{t,p}(\RN)} \leq C \|f\|_{L^m(\Omega)}.
$
\end{itemize}
Here, $C > 0$ are constants depending only on $N$, $s$, $t$, $p$, $m$ and $\Omega$.
\end{corollary}

\begin{proof}
The result immediately follows from Theorem \ref{calderon-zygmund-very-integrable}  and the fact that, for $0 < t < 1$ and $1 < p < \infty$,
$$
\|u\|_{L^{t,p}(\RN)} := \|u\|_{L^p(\RN)} + \|(-\Delta)^{\frac{t}{2}}u\|_{L^p(\RN)},
$$
is an equivalent norm for $L^{t,p}(\RN)$.
\end{proof}

\begin{corollary} \label{cor2-high-integrability}
Let $s \in (0,1)$, $\gamma \in (0,2s)$ and $u:= \mathbb{G}_s[f]$ with $f \in L^{m}(\Omega)$ for some $m > \max\big\{\frac{N}{s},\frac{N}{2s-\gamma}\big\}$: \begin{itemize}
\item[i)] If $0 < \gamma \leq s$, then $u \in W^{\gamma,p}(\RN)$ for all $1 < p < \infty$ and $\|u\|_{W^{\gamma,p}(\RN)} \leq C\, \|f\|_{L^m(\Omega)}$.
\item[ii)] If $s < \gamma < \min\{1,2s\}$, then $u \in W^{\gamma,p}(\RN)$ for all $1 < p < \frac{1}{\gamma-s}$ and $\|u\|_{W^{\gamma,p}(\RN)} \leq C \|f\|_{L^m(\Omega)}$.
\end{itemize}
Here, $C > 0$ are positive constants depending only on $N$, $s$, $\gamma$, $m$, $p$ and $\Omega$.
\end{corollary}

\begin{proof}
In the case where $0 < \gamma < s$, the result immediately follows combining Corollary \ref{cor1-high-integrability}, $i)$, with the embedding $L^{\sigma+\epsilon,p}(\RN) \subset W^{\sigma,p}(\RN)$ for all $0 < \epsilon < \sigma < 1$ and all $1 < p < \infty$ (see \cite[Theorem 7.63, (g)]{adams}).

To deal with the case where $s \leq \gamma < \min\{1,2s\}$, let us define
$$
\overline{p}(s,\gamma) := \left\{
\begin{aligned}
& + \infty, \quad && \textup{ for }\gamma = s, \\
& \frac{1}{\gamma-s}, \quad && \textup{ for } \gamma > s,
\end{aligned}
\right.
$$
and fix an arbitrary $1 < p < \overline{p}(s,\gamma)$. Since $m > \frac{N}{2s-\gamma}$ and $1 < p < \overline{p}(s,\gamma)$, we can choose
$$
\epsilon \in \Big(0, \min\Big\{\min\{1,2s\}-\gamma, 2s-\gamma-\frac{N}{m}, \frac{1}{p}-(\gamma-s) \Big\} \Big)
$$
and, by Corollary \ref{cor1-high-integrability}, $ii)$ applied with $t = \gamma + \epsilon$, it follows that
\begin{equation} \label{cor2-high-regularity-conclussion}
\|u\|_{L^{\gamma+\epsilon,p}(\RN)} \leq C\, \|f\|_{L^m(\Omega)}.
\end{equation}
The result in this case follows combining \eqref{cor2-high-regularity-conclussion} with the fact that $L^{\gamma+\epsilon,p}(\RN) \subset W^{\gamma,p}(\RN)$ for all $0 < \epsilon < \gamma < 1$ and all $1 < p < +\infty$.
\end{proof}

\begin{proof}[Proof of Theorem \ref{calderon-zygmund-low-integrable}]
The proof is similar to the proof of Theorem \ref{calderon-zygmund-very-integrable} so we omit some details. First of all, let us fix $g_1$, $g_2$ and $g_3$ as in \eqref{g1}, \eqref{g2} and \eqref{g3} respectively. Also, arguing as \eqref{far2-2}, for all $1 \leq p \leq +\infty$, we easily get $C := C (N,s,t,m,p,\Omega) > 0$ such that
\begin{equation} \label{th-1-3-conclussion2}
\|(-\Delta)^{\frac{t}{2}}u\|_{L^p(B_R^c(0))} \leq C \|f\|_{L^m(\Omega)}.
\end{equation}

Now, we consider separately two cases.
\medbreak
\noindent \textbf{Case 1 :} $m_{\star}(s,t) \leq m < \frac{N}{s}$.
\smallbreak
On one hand, combining Lemma \ref{Stein}, b) and c) with \eqref{dist2} and H\"older inequality, we infer that, for all $1 \leq \alpha < \frac{mN}{N-m(2s-t)}$, there exists $C > 0$ (depending only on $N$, $s$, $t$, $m$, $\alpha$ and $\Omega$) such that
\begin{equation} \label{th1-3-1}
\big\|g_1 + |\log \delta| \,g_2\big\|_{L^{\alpha}(\Omega)} \leq C \, \|f\|_{L^m(\Omega)}.
\end{equation}
On the other hand, combining Lemma \ref{Stein}, b) and c) with \eqref{dist} and H\"older inequality, we get that, for all $1 \leq \beta < \frac{mN}{N-ms+mN(t-s)} = p^{\star}(m,s,t)$, there exists $C = C(N,s,t,m,\b, \O)> 0$ such that
\begin{equation} \label{th1-3-2}
\Big\| \frac{g_3}{\delta^{t-s}} \Big\|_{L^{\beta}(\Omega)} \leq C \|f\|_{L^m(\Omega)}.
\end{equation}
Combining \eqref{th1-3-1} and \eqref{th1-3-2} with Lemma \ref{declem}, we conclude that, for all $1 \leq p < p^{\star}(m,s,t)$, there exists $C = C(N,s,t,m,p, \O)>0$
such that
\begin{equation} \label{th-1-3-conclussion1}
\|(-\Delta)^{\frac{t}{2}}u\|_{L^p(\Omega)} \leq C\, \|f\|_{L^m(\Omega)}.
\end{equation}
Having at hand \eqref{th-1-3-conclussion2} and \eqref{th-1-3-conclussion1}, it just remains to deal with $x \in B_R(0) \setminus \Omega$. Since $|x-y| \geq \max\{\delta(y),\delta(x)\}$ for all $y \in \Omega$ and $x \in B_R(0) \setminus \Omega$, for almost every $x \in B_R(0) \setminus \Omega$, it follows that
\begin{equation} \label{th-1-3-5}
\big| (-\Delta)^{\frac{t}{2}} u(x) \big| \leq \int_{\Omega} \frac{|u(y)|}{|x-y|^{N+t}} dy \leq \frac{1}{\delta^{t-s+\epsilon}(x)} \int_{\Omega} \frac{|u(y)|}{\delta^s(y)} \frac{dy}{|x-y|^{N-\epsilon}}, \quad  \textup{ for any } \epsilon > 0.
\end{equation}
Moreover, by Lemma \ref{hardy-reg}, for all $1 \leq \gamma < \frac{mN}{N-ms}$, there exists $C =C(N,s,t,m,\gamma, \O)>$
such that
\begin{equation} \label{th-1-3-6}
\Big\| \frac{u}{\delta^s} \Big\|_{L^{\gamma}(\Omega)} \leq C \, \|f\|_{L^m(\Omega)}.
\end{equation}
Combining \eqref{th-1-3-5} and \eqref{th-1-3-6} with Lemma \ref{Stein}, b) and c), \eqref{dist} and H\"older inequality, we get that, for all
$$
0 < \epsilon < \frac{N}{m}-s \quad \textup{ and } \quad 1 \leq q < \frac{mN}{N-ms+mN(t-s)+ \epsilon m (N-1)},
$$
there exists $C =C(N,s,t,m,\epsilon, q, \O)> 0$
such that
\begin{equation} \label{th-1-3-aux-conclussion}
\|(-\Delta)^{\frac{t}{2}}u\|_{L^q(B_R(0) \setminus \Omega)} \leq C \|f\|_{L^m(\Omega)}.
\end{equation}
Having at hand \eqref{th-1-3-aux-conclussion} we consider $1 \leq m < \frac{N}{s}$ and $1 \leq p < p^{\star}(m,s,t)$ fixed but arbitrary and choose
$$
\epsilon \in \Big(0, \min\Big\{\frac{N}{m}-s, \frac{1}{m(N-1)} \Big( \frac{mN}{p} - \big(N-ms+mN(t-s) \big) \Big) \Big\} \Big).
$$
Then, by \eqref{th-1-3-aux-conclussion}, we know there exists $C = C(N,s,t,m,p, \O)> 0$
such that
\begin{equation} \label{th-1-3-conclussion3}
\|(-\Delta)^{\frac{t}{2}}u\|_{L^p(B_R(0) \setminus \Omega)} \leq C \|f\|_{L^m(\Omega)}.
\end{equation}
Gathering \eqref{th-1-3-conclussion2}, \eqref{th-1-3-conclussion1} and \eqref{th-1-3-conclussion3} the result follows in the case where $1 \leq m < \frac{N}{s}$.
\medbreak
\noindent \textbf{Case 2 :} $\frac{N}{s} \leq m < \frac{N}{2s-t}$.
\smallbreak
Let us point out that, in this case, it is implicitly assumed that $s < t < \min\{1,2s\}$. By Lemma \ref{Stein}, b) and c) combined with \eqref{dist2} and H\"older inequality, it follows that, for all $1 \leq \alpha < \frac{mN}{N-m(2s-t)}$, there exists
$C= C(N,s,t, \a,m,\O) > 0$
such that
\begin{equation} \label{th1-3-3}
\big\|g_1 + |\log \delta| \,g_2\big\|_{L^{\alpha}(\Omega)} \leq C \, \|f\|_{L^m(\Omega)}.
\end{equation}
Moreover, combining Lemma \ref{Stein}, d) with \eqref{dist} and H\"older inequality, we get that, for all $1 \leq \beta < \frac{1}{t-s} = p^{\star}(m,s,t)$, there exists $C= C(N,s,t,m,\b,\O) > 0$
such that
\begin{equation} \label{th1-3-4}
\Big\| \frac{g_3}{\delta^{t-s}} \Big\|_{L^{\beta}(\Omega)} \leq C \|f\|_{L^m(\Omega)}.
\end{equation}
Since $m \geq \frac{N}{s}$ implies $\frac{mN}{N-m(2s-t)} > p^{\star}(m,s,t)$, combining \eqref{th1-3-3} and \eqref{th1-3-4} with Lemma \ref{declem}, we conclude that, for all $1 \leq p < p^{\star}(m,s,t)$, there exists $C= C(N,s,\gamma,m,\O) > 0$
such that
\begin{equation} \label{th-1-3-conclussion4}
\|(-\Delta)^{\frac{t}{2}}u\|_{L^p(\Omega)} \leq C\, \|f\|_{L^m(\Omega)}.
\end{equation}
Having at hand \eqref{th-1-3-conclussion2} and \eqref{th-1-3-conclussion4}, it remains to deal with $x \in B_R(0) \setminus \Omega$. Recall that, for almost every $x \in B_R(0) \setminus \Omega$, it follows that
\begin{equation} \label{th-1-3-7}
\big| (-\Delta)^{\frac{t}{2}} u(x) \big| \leq \int_{\Omega}  \frac{|u(y)|}{|x-y|^{N+t}} dy \leq \frac{1}{\delta^{t-s+\epsilon}(x)} \int_{\Omega} \frac{|u(y)|}{\delta^s(y)} \frac{dy}{|x-y|^{N-\epsilon}}, \quad \textup{ for any } \epsilon > 0.
\end{equation}
Moreover, by Lemma \ref{hardy-reg}, we know that, for all $1 \leq \gamma < + \infty$, there exists $C= C(N,s,\gamma,m,\O) > 0$
such that
\begin{equation} \label{th-1-3-8}
\Big\| \frac{u}{\delta^s} \Big\|_{L^{\gamma}(\Omega)} \leq C \, \|f\|_{L^m(\Omega)}.
\end{equation}
Combining \eqref{th-1-3-7} and \eqref{th-1-3-8} with Lemma \ref{Stein}, b) and c), \eqref{dist} and H\"older inequality, we get that, for all
$$
1 \leq \gamma < \infty, \quad 0 < \epsilon < \frac{N}{\gamma}, \quad \textup{ and } \quad 1 \leq q < \frac{N\gamma}{N+N\gamma(t-s) + \epsilon\gamma(N-1)},
$$
there exists $C= C(N,s,t,\gamma,\epsilon, q, m,\O) > 0$
such that
\begin{equation} \label{th-1-3-aux-conclussion-2}
\|(-\Delta)^{\frac{t}{2}}u\|_{L^q(B_R(0) \setminus \Omega)} \leq C \|f\|_{L^m(\Omega)}.
\end{equation}
Having at hand \eqref{th-1-3-aux-conclussion-2} we consider $\frac{N}{s} \leq m < \infty$ and $1 \leq p < p^{\star}(m,s,t)$ fixed but arbitrary and choose
$$
\gamma \in \Big( \frac{p}{1-p(t-s)}, + \infty \Big) \quad \textup{ and } \quad \epsilon \in \Big(0, \min \Big\{ \frac{N}{\gamma}, \frac{1}{\gamma(N-1)} \Big( \frac{N\gamma}{p}- \big(N+N\gamma(t-s) \big) \Big) \Big\} \Big).
$$
Then, by \eqref{th-1-3-aux-conclussion-2}, we know there exists $C= C(N,s,t,p, m,\O) > 0$
such that
\begin{equation} \label{th-1-3-conclussion6}
\|(-\Delta)^{\frac{t}{2}}u\|_{L^p(B_R(0) \setminus \Omega)} \leq C \|f\|_{L^m(\Omega)}.
\end{equation}
Gathering \eqref{th-1-3-conclussion2}, \eqref{th-1-3-conclussion4} and \eqref{th-1-3-conclussion6} the result follows also in this case.
\end{proof}

As a consequence of the above result, we deduce the following global regularity in Bessel potential spaces.

\begin{corollary} \label{cor1-low-integrability}
Let $s \in (0,1)$, $s \leq t < \min\{1,2s\}$ and $u:= \mathbb{G}_s[f]$ with $f \in L^{m}(\Omega)$ for some $m_{\star}(s,t) \leq m < \frac{N}{2s-t}$. Then, for all $1 < p < p^{\star}(m,s,t)$, there exists $C=C(N,s,t,m,p,\O) $  such that
$$
\|u\|_{L^{t,p}(\RN)} \leq C \|f\|_{L^m(\Omega)}.
$$
\end{corollary}

\begin{proof}
The result follows arguing again as in the proof of Corollary \ref{cor1-high-integrability} using now Theorem \ref{calderon-zygmund-low-integrable} instead of Theorem \ref{calderon-zygmund-very-integrable}.
\end{proof}

Next, we deal with the global Sobolev regularity.

\begin{corollary} \label{cor2-low-integrability}
Let $s \in (0,1)$, $s \leq \gamma < \min\{1,2s\}$ and $u:= \mathbb{G}_s[f]$ with $f \in L^{m}(\Omega)$ for some $m_{\star}(s,\gamma) \leq m < \frac{N}{2s-\gamma}$. Then, for all $1 < p < p^{\star}(m,s,\gamma)$, there exists $C=C(N,s,\gamma,m,p, \O) $ such that
$$
\|u\|_{W^{\gamma,p}(\RN)} \leq C \|f\|_{L^m(\Omega)}.
$$
\end{corollary}

\begin{proof}
Let $1 \leq m < \frac{N}{2s-\gamma}$ and $1 < p < p^{\star}(m,s,\gamma)$ fixed but arbitrary. Since $p < p^{\star}(m,s,\gamma)$, we can choose
$$
\epsilon \in \Big(0, \min\Big\{s,\min\{1,2s\}-\gamma,\frac{1}{p}-\frac{1}{p^{\star}(m,s,\gamma)} \Big\} \Big)
$$
and, by Corollary \ref{cor1-low-integrability}, applied with $t = \gamma+\epsilon$, it follows that
\begin{equation} \label{cor-Wsp-case1}
\|u\|_{L^{\gamma+\epsilon,p}(\RN)} \leq C \, \|f\|_{L^m(\Omega)}.
\end{equation}
The result follows combining \eqref{cor-Wsp-case1} with the fact that $L^{\gamma+\epsilon,p}(\RN) \subset W^{\gamma,p}(\RN)$ for all $0 < \epsilon < \gamma < 1$ and all $1 < p< +\infty$ (see e.g. \cite[Theorem 7.63, (g)]{adams}).
\end{proof}

\begin{proof} [Proof of Corollary \ref{cor-Wsp}]
It immediately follows from the combination of Corollaries \ref{cor2-high-integrability} and \ref{cor2-low-integrability}.
\end{proof}

Next, we prove the global weighted estimates.

\begin{theorem} \label{prop-distance-weight}
Let $s \in (0,1)$, $s \leq t < \min\{1,2s\}$, $u:= \mathbb{G}_s[f]$ with $f \in L^m(\Omega)$ for some $m \geq 1$ and
\begin{itemize}
\item[$\bullet$] $p = \infty$, if $m > \frac{N}{2s-t}$,
\item[$\bullet$] $1 \leq p < \infty$, if $m = \frac{N}{2s-t}$,
\item[$\bullet$] $1 \leq p < \frac{mN}{N-m(2s-t)}$, if $1 \leq m < \frac{N}{2s-t}$.
\end{itemize}
Then:
\begin{itemize}
\item[i)] If $t = s$, there exists $C > 0$ such that
$$
\big\| |\log \delta|^{-1}\,  (-\Delta)^{\frac{s}{2}} u \big\|_{L^p(\Omega)} \leq C \| f\|_{L^m(\Omega)}.
$$
\item[ii)] If $s < t < \min\{1,2s\}$, there exist $C > 0$ such that
$$
\big\|\delta^{t-s} (-\Delta)^{\frac{t}{2}} u \big\|_{L^p(\Omega)} \leq C \| f\|_{L^m(\Omega)}.
$$
\end{itemize}
Here, $C > 0$ are constants depending only on $N$, $s$, $t$, $p$, $m$ and $\Omega$.
\end{theorem}

\begin{proof}
We only detail the proof in the case where $t > s$, since the other estimates follows arguing on the same way. By Lemma \ref{declem} we know there exists $C= C(N,s,t, \O) > 0$
such that
$$
\delta^{t-s}(x) \big|(-\Delta)^{\frac{t}{2}} u(x) \big| \leq C \big[ g_1(x) + g_2(x) + g_3(x) \big], \quad \textup{ for a.e. }x \in \Omega.
$$
Here, $g_1$, $g_2$ and $g_3$ are as in \eqref{g1}, \eqref{g2} and \eqref{g3} respectively. Having at hand the above inequality, the result immediately follows from Lemma \ref{Stein}.
\end{proof}

\begin{proof}[Proof of Theorem \ref{weighted-distance-intro}]
It is a particular case of Theorem \ref{prop-distance-weight}.
\end{proof}

\begin{remark} The technique developed in this section allows to prove also estimates in $\O^c$. Indeed, suppose that the assumptions of Theorem \ref{prop-distance-weight} hold true and let $R$ be as in \eqref{erre}. As a byproduct of the proofs of Theorems \ref{calderon-zygmund-very-integrable} and \ref{calderon-zygmund-low-integrable}, we infer that, for $1 \leq p \leq \infty$, there exists $C > 0$ such that
$$
\big\|(-\Delta)^{\frac{t}{2}}u\big\|_{L^{p}(B_R^c(0))} \leq C\, \|f\|_{L^m(\Omega)}\,,
$$
On the other hand, taking into account \eqref{final---}, \eqref{out-weight-high-integrability}, \eqref{th-1-3-5} and \eqref{th-1-3-7} and arguing as in the proofs of Theorems \ref{calderon-zygmund-very-integrable} and \ref{calderon-zygmund-low-integrable}, one can infer that, for
\begin{itemize}
\item[$\bullet$] $1 \leq p \leq \infty$, if $m > \frac{N}{s}$,
\item[$\bullet$] $1 \leq p < \infty$, if $m = \frac{N}{s}$,
\item[$\bullet$] $1 \leq p < \frac{mN}{N-ms}$, if $1 \leq m < \frac{N}{s}$,
\end{itemize}
there exists $C > 0$ such that
$$
\big\|\, |\log \delta\,|^{-1} (-\Delta)^{\frac{s}{2}}u \big\|_{L^p(B_R(0)) \setminus \Omega)} +  (t-s)\big\|\,\delta^{t-s} (-\Delta)^{\frac{t}{2}} u \big\|_{L^p(B_R(0) \setminus \Omega)} \leq C \| f\|_{L^m(\Omega)}.
$$
\end{remark}

\bigbreak
Finally, let us provide the precise statement of the Calder\'on--Zygmund type regularity results for the \textit{Riesz fractional gradient}.

\begin{theorem} \label{CZ-gradient1}
Let $s \in (0,1)$, $s \leq t < \min\{1,2s\}$ and $u:= \mathbb{G}_s[f]$ with $f \in L^{m}(\Omega)$ for some $m > \frac{N}{2s-t}$. Then: \smallbreak
\begin{itemize}
\item[i)] For all $1 \leq p < \infty$, there exists $C > 0$ such that
$$
\|\nabla^s u\|_{L^p(\RN)} \leq C \|f\|_{L^{m}(\Omega)}.
$$
\item[ii)] For all $1 \leq p < \frac{1}{t-s}$, there exists $C > 0$ such that
$$
\|\nabla^t u\|_{L^p(\RN)} \leq C \|f\|_{L^{m}(\Omega)}.
$$
\end{itemize}
Here, $C > 0$ are constants depending only on $N$, $s$, $t$, $p$, $m$ and $\Omega$.
\end{theorem}

\begin{theorem} \label{CZ-gradient2}
Let $s \in (0,1)$, $s \leq t < \min\{1,2s\}$ and $u:= \mathbb{G}_s[f]$ with $f \in L^m(\Omega)$ for some $m_{\star}(s,t) \leq m < \frac{N}{2s-t}$. Then, for all $1 \leq p < p^{\star}(m,s,t)$, there exists $C > 0$ (depending only on $N$, $s$, $t$, $p$, $m$ and $\Omega$) such that
$$
\|\nabla^t u\|_{L^p(\RN)} \leq C\, \|f\|_{L^m(\Omega)}\,.
$$
\end{theorem}

\begin{remark}
Let us stress that the main results of our paper can be extended to more general operators than the $s$-Laplacian. The core of our proofs rely on the representation formula for the solution (see \eqref{representation-formula}) and on the estimates for the associated Green function and its gradient (see Lemma \ref{Green-estimates-lemma}). It these two ingredients are available for more general operators, our main regularity results can be extended to those operators. 

For instance, if we consider the operator $\mathcal{L}_s u :=(-\D)^s u+ B(x) \cdot \n u$ with $s \in (1/2,1)$ and $B \in (L^p(\Omega))^N$ for some $p > N/(2s-1)$, combining \cite[Theorem 1]{BJ} and \cite[Corollary 1.8]{KUL} with Lemma \ref{Green-estimates-lemma}, we infer that the Green function associated to $\mathcal{L}_s$ satisfies \eqref{green100}--\eqref{green200}. 
Hence, the unique solution to 
$$
 \left\{
\begin{aligned}
\, \mathcal{L}_s u &= f, && \text{ in }\Omega , \\
u &= 0, &&\hbox{  in } \mathbb{R}^N\setminus\Omega,\\
\end{aligned}
\right.
$$ 
with $f \in L^{m} (\Omega)$, $m\geq1$,  satisfies the estimates of Theorem \ref{calderon-zygmund-very-integrable} and  Theorems \ref{prop-distance-weight}--\ref{CZ-gradient2}.  
\end{remark}

\appendix

\section{Proof of Proposition \ref{integrability-t-Laplacian-intro}}

\begin{proof}[Proof of Proposition \ref{integrability-t-Laplacian-intro}] By \cite[Theorem 1]{Dyda} we know that, for all $0 < t < 2$, 
\begin{equation*}
(-\Delta)^{\frac{t}{2}}\bu(x) = - \frac{2^{t-2s} \Gamma \big( \frac{N+t}{2}  \big)}{\Gamma \big(\frac{N}{2}+s \big) \Gamma \big( s+1-\frac{t}{2} \big)} \, \,{}_2 F_1 \left( \frac{t+N}{2}, -s + \frac{t}{2}; \frac{N}{2}; |x|^2 \right), \quad \mbox{ for all } x \in B_1(0)\,.
\end{equation*}
Here, ${}_2 F_1$ denotes the Gauss' hypergeometric function.  
Moreover, it is well known (see e.g. \cite[Theorem 2.1.3]{AAR1999}) that, for all $s <t < \min\{2s,1\}$,
\begin{equation*} 
\lim_{|x| \to 1^{-}} \frac{{}_2 F_1 \left( \frac{t+N}{2}, - s + \frac{t}{2}; \frac{N}{2}; |x|^2 \right)}{ (1-|x|^2)^{s-t}} = \frac{\Gamma\big( \frac{N}{2} \big) \Gamma\big(t-s\big)}{ \Gamma \big( \frac{N+t}{2} \big) \Gamma \big( -s+\frac{t}{2} \big) }, 
\end{equation*}
and that, for $s = t$, 
\begin{equation*} 
\lim_{|x| \to 1^{-}} \frac{{}_2 F_1 \left( \frac{t+N}{2}, - s + \frac{t}{2}; \frac{N}{2}; |x|^2 \right)}{ - \ln(1-|x|^2)} = \frac{\Gamma\big( \frac{N}{2} \big)}{ \Gamma \big( \frac{N+s}{2} \big) \Gamma \big(-\frac{s}{2} \big) }, 
\end{equation*}
Using this information and direct (elementary) computations, we deduce that Proposition \ref{integrability-t-Laplacian-intro} holds true.  \end{proof}

As a byproduct of Proposition \ref{integrability-t-Laplacian-intro}, the following regularity in  Bessel, Sobolev and Besov potential spaces can be obtained. 

\begin{corollary} 
$ $
\begin{enumerate}
\item[i)] If $0 < t < s$, then $\bu \in W^{t,p}(\RN) \cap \bu \in L^{t,p}(\RN)$ for all $1 < p < \infty$. 

\item[ii)] If $t = s$, then:
\begin{itemize}
\item[a)] $\bu \in L^{s,p}(\RN)$ for all $1 < p < \infty$.
\item[b)] $\bu \in B^{s}_{p,2}(\RN)$ for all $1 < p < 2$,  while 
 $\bu \in W^{s,p}(\RN)$ for all $2 \leq p < \infty$.
\end{itemize}
\item[iii)] If $s < t < \min\{1,2s\}$, then:
\begin{itemize}
\item[a)]  $\bu \in L^{t,p}(\RN)$ for all $1 < p < \frac{1}{t-s}$ and $\bu \not \in L^{t,q}(\RN)$ for $ \frac{1}{t-s} \leq q < \infty$.
\item[b)] $\bu \in B^{t}_{p,2}(\RN)$ for all $1 < p < 2$, while
 $\bu \in W^{t,p}(\RN)$ for all $ 2 \leq p < \frac{1}{t-s}$.

\end{itemize}
\end{enumerate}
\end{corollary}

\vspace{0.2cm}
\bibliographystyle{plain}
\bibliography{Bibliography}
\vspace{0.2cm}

\end{document}